\documentclass[11pt, reqno]{amsart}

\newcommand{\mysection}[1]{\section{#1}
      \setcounter{equation}{0}}

\renewcommand\){{\rm)}}

\newcommand\loc{\textnormal{loc}}

\newcommand{\sign}{\text{\rm\,sign}\,}

\newtheorem{theorem}{Theorem}[section]
\newtheorem{lemma}[theorem]{Lemma}

\theoremstyle{definition}
\newtheorem{assumption}{Assumption}[section]

\theoremstyle{remark}
\newtheorem{remark}{Remark}[section]

\newtheorem{example}{Example}[section]

\newcommand\cbrk{\text{$]$\kern-.15em$]$}}
\newcommand\opar{\text{\raise.2ex\hbox{${\scriptstyle | }$}\kern-.34em$($} }
\newcommand{\tr}{\text{\rm tr}\,}

  \makeatletter
 \def\dashint{%
 \operatorname%
 {\,\,\text{\bf--}\kern-.98em\DOTSI\intop\ilimits@\!\!}}
 \makeatother
 \newcommand{\WO}{\overset{\rm o}{ W}\,\!}
 
\newcommand\bR{\mathbb{R}}

\newcommand\bL{\mathbb{L}}

\newcommand\bS{\mathbb{S}}

\newcommand\bZ{\mathbb{Z}}

\newcommand\cP{\mathcal{P}}

\newcommand\cH{\mathcal{H}}

\newcommand{\Span}{{\rm Span}\,}
\newcommand{\diam}{{\rm diam}}
\newcommand\dist{{\rm dist}\,}

\begin{document}

\title[Fully nonlinear
elliptic  
equations]
{On the existence of $W^{2}_{p}$ 
solutions for fully nonlinear elliptic
equations under either relaxed 
or no convexity assumptions}

\author{N.V. Krylov}
\thanks{The  author was partially supported by
 NSF Grant DMS-1160569  and
by a grant 
from the Simons Foundation (\#330456 to Nicolai Krylov)}
\email{krylov@math.umn.edu}
\address{127 Vincent Hall, University of Minnesota,
 Minneapolis, MN, 55455}



 \begin{abstract}
We establish the existence  of solutions
of fully nonlinear elliptic second-order equations
like $H(v,Dv,D^{2}v,x)=0$
in smooth domains without requiring $H$
to be convex or concave with respect to the second-order derivatives.
Apart from ellipticity nothing is required of $H$ at points at which
 $|D^{2}v|\leq K$,
where $K$ is any given constant. For large $|D^{2}v|$ some kind
of relaxed convexity assumption with respect to $D^{2}v$ mixed with a VMO condition
with respect to $x$ are still imposed. The solutions are sought
in Sobolev classes. We also establish the solvability
without almost any conditions on $H$, apart from ellipticity,
but of a ``cut-off'' version of the equation $H(v,Dv,D^{2}v,x)=0$.
\end{abstract}

\maketitle

\mysection{Introduction and main results}

The first object, we deal with  in this paper, is
the equation
\begin{equation}
                                                \label{7.29.1}
H[v](x):= H(v( x),D v( x),D^{2}v( x), x)=0
\end{equation}
considered in subdomains of $\bR^{d } $, 
where $H(u,x)$ is a real-valued
function defined for
$$
u=(u',u''),\quad u'=(u'_{0},u'_{1},...,u'_{d}),\quad u''\in\bS,
\quad x\in\bR^{d}, 
$$
where $\bS$ is the set of symmetric $d\times d$ 
matrices. Let $\Omega$ be an open bounded subset of $\bR^{d}$,
fix
  $p>d$ and   functions
  $\bar G,f\in L_{p}(\Omega)$, $\bar G\geq0$.
One of our
 main results implies that, for $d=3$ and  $\Omega\in C^{2}$,  the equation
($a\wedge b=\min(a,b)$)
$$
H(D^{2}u,x):=\bar G(x)\wedge|D_{12}u|
+\bar G(x)\wedge|D_{23}u|+\bar G(x)\wedge|D_{31}u|
$$
\begin{equation}
                                                      \label{3.4.1}
+3\Delta u-f(x)=0
\end{equation}
in $\Omega$ with zero boundary condition has a unique
solution $u\in W^{2}_{p}(\Omega)$.
Recall that $W^{2}_p(\Omega)$ denotes the set of
functions $v$ defined in
$\Omega$ such that $v$, $Dv
 $, and $D^2v $ 
are in
$L_p(\Omega)$. Observe that $H$ in \eqref{3.4.1}
is neither convex nor concave with respect to $D^{2}u$.
So far, there are only two approaches to such equations:
the theory of ($L_{p}$) viscosity solutions
and the theory of stochastic differential games,
provided $H$ has a somewhat special form. The past experience
shows that
it is hard
to expect getting sharp quantitative results using probability theory.
On the other hand, the theory of viscosity solutions indeed produced
some remarkable quantitative results.
However, to the best of the author's knowledge the result
stated above about \eqref{3.4.1}
is either very hard to obtain by using the theory
of ($L_{p}$) viscosity solutions or is just beyond it,
at least at the current stage.
It seems that the best information, that theory provides
at the moment,
is the existence of the maximal and minimal
$L_{p}$-viscosity solution (see \cite{JS_05}),
no uniqueness of $L_{p}$-viscosity solutions
can be inferred  for \eqref{3.4.1} and no regularity
apart from the classical $C^{\alpha}$-regularity (see
\cite{CKLS99}).

 Fix some constants $K_{0},K_{F}\in[0,\infty)$,
$\delta\in(0,1]$,
and an increasing
 continuous function $\omega(r)$, $r\geq0$,
such that $\omega(0)=0$.

 \begin{assumption}
                                    \label{assumption 10.30.1}
 
The function   $H(u, x)$
is measurable.
  
 \end{assumption}

Next, we assume that there are two
functions $F(u,x)=F(u'',x)$ and $G(u,x)$ such that
$$
H =F +G .
$$
 
\begin{example}
                                          \label{example 11.29.1}
One can take
$F(u'',x)=H(0,u'',x)$ and $G=H-F$. 
Since we will require later that $F(0,x)=0$, one can then
take $F(u'',x)=H(0,u'',x)-H(0,x)$ and $G=H-F$. However, we are
not bound by these options.
\end{example}

 The following assumptions contain   
parameters $\hat\theta,\theta\in(0,1]$
which are specified later in our results. 

\begin{assumption}
                                    \label{assumption 10.5.1}

 For $u''\in\bS,u'\in\bR^{d+1}$, and $x\in\bR^{d}$ we have
$$
|G(u',u'',x)|\leq  \hat\theta |u''|+ K_{0}|u'|+\bar{G}(x).
$$

\end{assumption}
   
 Set
$$
B_{r}(x)=\{y\in\bR^{d}:|x-y|<r\},\quad B_{r}=B_{r}(0),
$$
and for Borel $\Gamma\subset \bR^{d}$ denote by $|\Gamma|$
the volume of $\Gamma$.
Also set
$$
\bS_{\delta}=\{a\in\bS:\delta|\xi|^{2}\leq a_{ij}\xi_{i}\xi_{j}
\leq\delta^{-1}|\xi|^{2},\quad \forall \xi\in\bR^{d}\}.
$$
 
Recall that   Lipschitz
continuous functions are almost everywhere differentiable.

\begin{assumption}
                                \label{assumption 2.13.4}
(i\,)  The function $F$ is Lipschitz continuous with respect to $u''$
with Lipschitz constant $K_{F}$, measurable
with respect to $x$, and 
$$
F(0,x)\equiv0.
$$

Moreover, there exist  $R_0\in(0,1]$ and $t_{0}\in[0,\infty)$
such that for any
$r\in (0, R_0]$  and  $z\in \Omega$
one can find a {\em convex\/} function $\bar{F} (u'' )=
\bar{F}_{z,r} (u'' )$ (independent
of $x $)  such that

(ii\,) We have $\bar{F}(0)=0$ and at all points of differentiability
of $\bar{F}$ we have $(\bar{F}_{u''_{ij}})\in\bS_{\delta}$;
 
(iii\,) 
For any $u''\in\bS$ with $|u''|=1$,   we have
\begin{equation}
                                                \label{2.13.5}
\int_{ B_{r}(z)\cap\Omega}\sup_{t>t_{0}}t^{-1}
|F(tu'' ,x)-\bar{F}(tu'')| \,dx\leq \theta
| B_{r}(z)\cap\Omega|.
\end{equation}
 
\end{assumption}

\begin{assumption}
                                \label{assumption 3.10.4}
We have $\Omega\in C^{2}$, and we are given
a function $g\in W^{ 2}_{p }(\Omega) $.
\end{assumption}

We only consider $p>d$ and then by embedding theorems
the functions of class $W^{ 2}_{p }(\Omega)$
admit  modifications that are continuous in $\bar \Omega$
along with their first derivatives.
We will always have in mind such modifications.

Here are our first two main results.
The first one is an  a priori estimate.
Observe that   in Theorem \ref{theorem 2.20.1}
we do not even assume that
equation  \eqref{7.29.1} is elliptic.

\begin{theorem}
                                            \label{theorem 2.20.1}
Under the above assumptions
there exists a constant $\theta\in (0,1]$, depending
only on $d$, $p$, $\delta$, and  $  \Omega$,
  and a constant $\hat\theta\in(0,1]$,
depending only on  $K_{0}$, $K_{F}$, $d$, $p$, $\delta$,   
$R_{0}$, and
  $  \Omega $, 
such that, if  
Assumptions \ref{assumption 2.13.4} and \ref{assumption 10.5.1}
 are satisfied with these
 $\theta$ and $\hat\theta$, respectively,   then for any 
 $u\in W^{ 2}_{d }(\Omega) $,
 that satisfies \eqref{7.29.1} in $\Omega$ \(a.e.\)
and equals $g$ on $\partial\Omega$, we have
\begin{equation}
                                \label{eq16.01}
\|u\|_{W^{2}_p(\Omega )} \le N \big(
\| \bar G\|_{L_p(\Omega)}+\|g\|_{W^{2}_p(\Omega )}+
 \|u\|_{L_p(\Omega)}\big)
 +N t_{0} , 
\end{equation}
where $N$ depends only on
 $K_{0}$, $K_{F}$, $d$, $p$, $\delta$,   
$R_{0}$, and
  $  \Omega $.
\end{theorem}

In Section \ref{section 3.14.1}
we will see that Theorem \ref{theorem 2.20.1}
is an easy generalization of the estimate
in Theorem 2.1 of \cite{Kr13}. However, in \cite{Kr13}
the ellipticity of equation \eqref{7.29.1} is assumed from
the start. It is also supposed and used that Assumption \ref{assumption 3.11.2}
below is satisfied with $\omega(t)=K_{0}t$, and thus $N$
in the estimate corresponding to \eqref{eq16.01}  in
Theorem 2.1 of \cite{Kr13}   depends on the constant
of Lipschitz continuity
 of $H$ with respect to $u'$. We will see that, actually,
it is independent.

 To have the solvability we need ellipticity
and more regularity of $H$. 
 
\begin{assumption}
                                 \label{assumption 3.11.2}
\(i\,\) The function   $H(u, x)$
is  
  Lipschitz continuous with respect to $u''$, and at 
all points of differentiability
of $H$ with respect to $u''$ we have
$D_{u''}H\in \bS_{\delta}$;

\(ii\,\)
 The function   $H(u, x)$
 is nonincreasing with respect to $u'_{0}$ and
$$
|H(u',u'', x)-H(v',u'', x)|\leq\omega(|u'-v'|)
$$
for all $u,v $, and $x$.
\end{assumption}

\begin{theorem}
                                    \label{theorem 2.20.2}
 
Suppose that the above assumptions are satisfies and, moreover,
Assumptions \ref{assumption 10.5.1} and  \ref{assumption 2.13.4} 
are satisfied with  
  $\hat\theta$ and   $\theta$ 
from Theorem \ref{theorem 2.20.1}.

Then
for any $g\in W^{ 2}_{p }(\Omega) $ there exists  
$u\in W^{ 2}_{p }(\Omega) $
 satisfying \eqref{7.29.1} in $\Omega$ \(a.e.\)  and such
that $u=g $ on $\partial\Omega$.

\end{theorem}

Theorem \ref{theorem 2.20.2}
is a generalization of the existence part of
Theorem 2.1 of \cite{Kr13}. Again the whole point
is that here we do not assume that $H$
is Lipschitz in $u'$, which makes some arguments
more involved. In particular, in \cite{Kr13}
the existence is first proved when $g=0$
and then the unknown $u$ is replaced with $u-g$.
In our present case this would cause the last
$N$ in \eqref{eq16.01} to also depend on $\|g\|_{W^{2}_p(\Omega )}$.
 Then, the results of 
\cite{Kr13} are based on Theorem 1.1 of \cite{Kr12}
in which the assumption of the Lipschitz
continuity of $H$ with respect to $u'$
was indispensable. We abandon this assumption
and then we have to show that the corresponding
counterpart of Theorem 1.1 of \cite{Kr12},
which is our third main result, Theorem \ref{theorem 8.29.1},
still holds.

These results are natural continuation
of the results in \cite{Kr13}.
Some discussion and example of applications of
such results can be found  there.
In particular, an application to Isaacs equations
is given.

\begin{remark}
                                \label{remark 3.4.1}
In case of equation \eqref{3.4.1}, to prove what is claimed after
it, it suffices to take $F(D^{2}u,x)=\bar F(D^{2}u)=
3\Delta u$,
and then all our assumptions are obviously satisfied.
In particular, Assumption \ref{assumption 3.11.2} (i)
is satisfied since $H_{u''_{ii}}=3$, $i=1,2,3$, $|H_{u''_{12}}|=
|H_{u''_{21}}|,
|H_{u''_{23}}|=|H_{u''_{32}}|,|H_{u''_{31}}|
=|H_{u''_{13}}|\leq 1$.
\end{remark}

 In the literature,   interior $ W^2_p, p>d$, {\em a priori\/}
estimates for a class of
fully nonlinear uniformly elliptic equations in $\bR^{d}$ of the form
\eqref{7.29.1} in terms of viscosity solutions were first obtained by Caffarelli in \cite{Caf89} (1989) 
(see also
\cite{CC95} (1995)).
 Adapting his technique, similar interior 
a priori estimates were proved by Wang
\cite{Wa92} (1992) for parabolic equations. In the same paper, a boundary
estimate is stated but without   proof; see Theorem 5.8 there.  By
exploiting a weak reverse H\"older's inequality, the result of
\cite{Caf89} was sharpened by Escauriaza in \cite{Es93} (1993), 
who obtained the
interior $ W^2_p$-estimate for the same equations allowing
$ p>d-\varepsilon$, with a small constant
$\varepsilon$ depending only on the ellipticity constant and $d$. 

The above cited works    are quite remarkable 
in one   respect--they do not suppose that $H$ is convex
 or concave in $ D^{2}u$.
But  they only show that 
to prove a priori estimates
it suffices to prove 
the interior $  C^{2}$--estimates for
``harmonic'' functions. However, up to now,
these estimates  are only known under convexity assumptions.
 
Generally, having a priori estimates does not guarantee
that there are existence results. Say, if $d=1$
 and equation is
$$
  u''-(2u'')I_{|u''|  \leq   5}-u=1\quad\text{in}\quad (-1,1).
$$
a priori $W^{2}_{p}$-estimates are trivial since $|(2u'')I_{|u''|  \leq   5}|
\leq10$. However, the solvability is quite questionable.

Also obtaining boundary $ W^{2}_{p}$ estimate by using the theory
of viscosity solutions
 turned out to be  extremely challenging and only
in 2009, twenty years after the work of Caffarelli,
 Winter \cite{Wi09}   proved the solvability
 in $ W^{2}_{p}(\Omega)$ of equation \eqref{7.29.1} 
 with Dirichlet boundary condition in $ \Omega\in C^{1,1}$. 

It is
also worth noting that a solvability theorem in the space
$  W^{1,2}_{p,\text{ loc}}(Q)\cap C(\bar Q)$ can be found in
M. G. Crandall, M. Kocan, A. \'Swi{\c e}ch \cite{CKS00} (2000) for
the boundary-value problem for fully nonlinear parabolic equations. 
The above mentioned existence results of \cite{CKS00} and \cite{Wi09}
are proved under the assumption that $ H$ is convex in $ D^{2}v$
and in all papers mentioned above
 a small oscillation assumption in the integral sense is
imposed on the operators; see below. 
The above cited works  are performed in the framework of
viscosity solutions.

  Caffarelli and Cabr\'e \cite{CC95} 
consider  equations
$$
 H(D^{2}u,x){=}f(x) 
$$
with $ H(0,x){=}0$.
They, and a very many other authors after them, introduced
$$ 
\beta(x,x_{0}){=}\sup_{u''\in\bS}\frac{|H(u'',x)-H(u'',x_{0})|}{|u''|}.
$$
and require that for all $ B_{r}(x_{0})\subset B_{1}$
\begin{equation} 
                                                 \label{10.24.1}
\dashint_{B_{r}(x_{0})}\beta^{d}(x,x_{0})\,dx
\,{\leq}\, \theta^{d}(p),
\end{equation}
where $ p>d$ and $ \theta(p)$ is small enough. Then they prove,
 under an additional assumption on $ H(u'',0)$, that
the $ W^{2}_{p}(B_{1/2})$-norm of
the solution can be controlled.
It is easily shown that
$ \theta(p)\downarrow 0$ as $ p\to{\infty}$ implies that $ H(u'',x)$
is just a continuous function of $ x$.

Condition \eqref{10.24.1} might look like 
our condition \eqref{2.13.5}. Therefore, it is important
to emphasize that they are quite different.
For instance, in the case of equation \eqref{3.4.1}, if $\bar G\leq 1$,
condition \eqref{10.24.1} implies that 
 $$
  \bigg|\dashint_{B_{r}(x_{0})}G(x)\,dx-G(x_{0})\bigg|\leq\theta(p),
$$
where the integral is a uniformly continuous function of $x_{0}$.
According to Remark \ref{remark 3.4.1} we do not need any conditions
on $\bar G$ apart from $\bar G\in L_{p}(\Omega)$.

Probably, even better illustration of the difference gives
the example of {\em linear\/} equations when $H[u]=a^{ij}D_{ij}u$.
Then condition \eqref{10.24.1} is equivalent to the fact that
$a^{ij}(x)$ are uniformly sufficiently close to uniformly continuous
ones, and our condition is satisfied if, say $a^{ij}\in VMO$.

Also the additional assumption in 
\cite{Caf89}, 
\cite{CC95} 
on $H(u'',x)$ alluded to above,
in the case of equation \eqref{3.4.1}, includes the requirement 
that, for $f=0$, \eqref{3.4.1} with $\bar G(0)$ in place of
$\bar G(x)$ admit a solution of class $C^{2}(B_{1/2})$
for any continuous boundary data. So far, we have no idea whether
  this happens indeed and this is a very challenging problem.

 Observe that 
 in \cite{CKS00} and \cite{Wi09} the assumption that $H$ is convex
in $u''$ is needed because for such equations with additional
regularity assumption the solvability in $C^{2+\alpha}$ is known
and one can approximate  
the solutions  with more regular ones in a usual way.

Our approach to construct the approximations is different and
is based on  our third main result before which we introduce
some notation and assumptions.

Here is Theorem  3.1 of \cite{Kr11}.

\begin{theorem}
                                               \label{theorem 9.26.3}
There exists an integer $m=m(\delta,d)>d$, 
a set 
\begin{equation}
                                              \label{6.1.2}
\Lambda=\{ l_{\pm1},...,l_{\pm m}\}
\subset \bar B_{1}\subset\bR^{d},
\end{equation}
 and 
 there exists 
a constant  
$$
\hat{\delta}=\hat{\delta}(\delta,d ) 
\in(0,\delta/4]
$$
such that the coordinates of $l_{  k}$ are rational numbers and

\(a\,\) We have $l_{-k}=-l_{k}$, $k=1,...,m$,
$ e_{i}, (1/2)(e_{i}+e_{j})\in\Lambda$, $i\ne j$.
  $i,j=1,...,d$, where $e_{1},...,e_{d}$ is the 
 standard orthonormal basis of $\bR^{d}$;

\(b\,\) There   exist
real-analytic functions $\lambda_{\pm 1}(a),...,
\lambda_{\pm m}(a)$ on $\bS_{\delta/4}$ such that
for any $a\in\bS_{\delta/4}$
\begin{equation}
                                         \label{4.8.01}
a=\sum_{|k|=1}^{m}\lambda_{k}(a)l_{k}l_{k}^{*},
\quad \hat{\delta}^{-1}\geq\lambda_{k}(a)\geq\hat{\delta}
,\quad \forall k.
\end{equation}
\end{theorem}

For $z''=(z''_{\pm 1},...,z''_{\pm m})$ introduce
\begin{equation}
                                                       \label{6.1.01}
\cP(z'')  =\max_{\substack{\hat{\delta}/2\leq a_{k}\leq
2\hat{\delta}^{-1} \\    |k|=1,...,m} }
 \sum_{|k|=1}^{m} a_{k}
z''_{k} ,
\end{equation}
and for $ u'' \in \bS$ define
\begin{equation}
                                          \label{3.28.1}
P( u'')=\cP( \langle u''l_{\pm 1},l_{\pm 1}\rangle,...,
\langle u''l_{\pm m},l_{\pm m}\rangle),
\end{equation}
where $\langle\cdot\,,\cdot\rangle$ is the scalar product in $\bR^{d}$.

\begin{remark}
                                                \label{remark 11.1.2}
Observe (cf.~Section \ref{section 3.19.1}) 
that  $P( u'')$ is Lipschitz continuous
and at all points of its differentiability, for any $\xi\in\bR^{d}$, we have
$$
\min_{\substack{\hat{\delta}/2\leq a_{k}\leq
2\hat{\delta}^{-1} \\    |k|=1,...,m} }
 \sum_{|k|=1}^{m} a_{k}\langle l_{k},\xi\rangle^{2}
\leq\langle \xi,P_{u''}(u'')\xi\rangle\leq
\max_{\substack{\hat{\delta}/2\leq a_{k}\leq
2\hat{\delta}^{-1} \\    |k|=1,...,m} }
 \sum_{|k|=1}^{m} a_{k}\langle l_{k},\xi\rangle^{2}.
$$
Since, $e_{i}\in\Lambda\subset\bar B_{1}$, we see that
there is a $\bar \delta=\bar\delta(d,\delta)\in(0,\delta]$
such that $P_{u''}(u'')\in\bS_{\bar\delta}$.
\end{remark}

For smooth enough functions $u(x)$ introduce
$$
P[u](x)
=P( D^{2}u(x)).
$$
\begin{assumption}
                                \label{assumption 8.21.1}

We are given a domain $\Omega
\subset\bR^{d}$ which is bounded and satisfies the exterior
ball condition. We are also
 given a function $g\in C^{1,1}(\bR^{d})$.
\end{assumption} 

 By the exterior ball condition we mean that
for any $x_{0}\in\partial\Omega$ there exists
a closed ball of radius $\rho(\Omega)$ (independent of $x_{0}$ and $\leq1$)
centered outside $\Omega$ and having $x_{0}$
as the only common point with $\bar \Omega$.

\begin{assumption}
                                    \label{assumption 8.30.1}
(i) 
Assumptions \ref{assumption 10.30.1} and  \ref{assumption 3.11.2} 
are satisfied.

(ii) The number
$$
\bar{H} :=\sup_{u',  x}\big(|H (u',0, x)|-K_{0}|u'|\big)
\quad(\geq0)
$$
is finite.

 \end{assumption}

Fix a constant $K\geq0$ and consider the equation
\begin{equation}
                                             \label{9.23.02}
 \max(H[v],P[v]-K)=0 \quad\text{in}\quad\Omega\quad\rm (a.e.)
\end{equation}
with boundary condition $v=g$ on $\partial\Omega$.
Here is the major result concerning \eqref{9.23.02}.
It generalizes Theorem 1.1 of \cite{Kr12} by relaxing the
requirement of the Lipschitz continuity of $H$
with respect to $u'$ to just continuity. As there,
we do not assume any regularity of the dependence
of $H$ on $x$ and still obtain solutions with locally bounded
second-order derivatives.
Set
$$
\rho_{\Omega}(x)=\dist(x,\bR^{d}\setminus\Omega) .
$$

\begin{theorem}
                                    \label{theorem 8.29.1}
Let Assumptions 
\ref{assumption 8.21.1} and \ref{assumption 8.30.1}
be satisfied. Then 

\(i\)
equation \eqref{9.23.02}
 with boundary condition $v=g$ on $\partial\Omega$
has a solution  $v\in
C(\bar{\Omega} )\cap C^{1,1}_{ \loc }(\Omega)$
such that, 
\begin{equation}
                                             \label{10.18.3}
|v|\leq N(\bar H+\|g\|_{ C (\Omega)}),
\end{equation}
in $\Omega$, where $N$ is a constant
depending only on $d,\delta,K_{0}$, and $\diam(\Omega)$, and,
for all
$i,j=1,...,d$, 
\begin{equation}
                                             \label{8.29.10}
|D_{i}v|, \rho_{\Omega} |D_{ij} v | \leq N(\bar{H}+K
+\|g\|_{ C^{1,1} (\bR^{d})})\quad\text{in} 
\quad \Omega  \quad (a.e.),
\end{equation}
where 
   $N$ is a constant depending only on $\Omega$,
  $K_{0}$,   $\Lambda$, and $\delta$
\(in particular, $N$ is independent
of $\omega$\);

\(ii\) for $p>d$ and any such solution
\begin{equation}
                                                \label{1.13.2}
\|v\|_{W^{2}_{p}(\Omega)}\leq
N_{p}(\bar{H}+K+\|g\|_{ C^{1,1} (\bR^{d})}),
\end{equation}
where $N_{p}$ is a constant depending only on $p$, $\Omega$,
  $K_{0}$,   $\Lambda$, and $\delta$.

\end{theorem}

Results of this kind, valid for $H$ with no concavity or convexity
assumptions, have independent interest and have already been
used in \cite{Kr_15} to show that the value functions
in stochastic differential games in $\Omega$ admit
approximations of order $1/K$ by functions whose second-order
derivatives by magnitude 
are of order $K$ in any compact subdomain of $\Omega$.
They were also used to show that solutions of Isaacs
equations are in $C^{1+\alpha}$ if the coefficients
are in VMO (see \cite{Kr_14}),
and that in a rather general case solutions of elliptic
and parabolic Isaacs equations admit unique viscosity
solutions that can be approximated, by
using finite-difference equations,
with algebraic rate of convergence with respect
to the mesh size (see \cite{Kr_14_1} and \cite{Kr_15_1}).

Here we use Theorem \ref{theorem 8.29.1} to prove
Theorem \ref{theorem 2.20.2}.

\begin{remark}
                                         \label{remark 3.10.1}
Generally, in the situations of Theorem
\ref{theorem 2.20.2} and \ref{theorem 8.29.1}
there is no uniqueness. For instance,
in the one-dimensional case each of
the  equations 
$$
 D^{2}u -\sqrt{12|Du |} =0,\quad
\max(D^{2}u -\sqrt{12|Du |}, 2(D^{2}u)^{+}-
(1/2)(D^{2}u)^{-} -7)=0
$$
 for $x\in(-1,1)$
with zero boundary data has two solutions:
one is identically equal to zero and the other one is
 $1-|x|^{3}$.

To guarantee uniqueness, it suffices to assume that
$\omega(t)=N_{0}t$, where $N_{0}$ is a constant. This follows from the fact that
if $u,v$ are solutions of, say \eqref{9.23.02}
with the properties described in Theorem \ref{theorem 8.29.1} (i),
then (cf.~Section \ref{section 3.19.1}) the function $w=u-v$
equals zero on $\partial\Omega$ and in $\Omega$ satisfies
the equation
$$
a^{ij}D_{ij}w+b^{i}D_{i}w-cw=0,
$$
where $(a^{ij})$ is an $\bS_{\bar\delta}$-valued function,
$|(b^{i})|\leq N_{0}$, $0\leq c\leq N_{0}$. Then the equality
$w=0$ follows from the Aleksandrov maximum principle.

\end{remark}

As we have pointed out already,
the above results are close to those   in \cite{Kr13}
and \cite{Kr12}, where the assumptions
are stronger than here. 
The major difference between the assumptions in \cite{Kr13}
and \cite{Kr12}
and in the present article is that we do not
assume that $H$ is Lipschitz continuous with respect
to $u'$ with Lipschitz constant independent of $u'',x$.
Thus, we considerably enlarge the number of equations
for which the existence of solutions is guaranteed.
 
Abandoning the Lipschitz continuity with respect to
$u'$ causes some serious complications in the proof
of our main start-up result: Theorem \ref{theorem 8.29.1}.
In \cite{Kr13} we  used an auxiliary cut-off equation
\begin{equation}
                                                 \label{3.4.4}
\max(H[v],P(u,Du,D^{2}u)-K)=0
\end{equation}
and used a result, proved in \cite{Kr12}, that   such equations
have solutions with locally bounded second-order derivatives. The reason 
why it happens is that, owing to the construction
of $P$, on the set where
\begin{equation}
                                                 \label{3.6.2}
H(v( x),D v( x),D^{2}v( x), x)\geq P(u,Du,D^{2}u)-K,
\end{equation}
we had $|v|,|Dv|,|D^{2}v|\leq N(K)$. On the complement
of this set
\begin{equation}
                                                 \label{3.6.1}
P(u,Du,D^{2}u)-K=0.
\end{equation}
In addition, the constructed $P$ in \cite{Kr12}
(and here) is convex in $v,Dv,D^{2}v$, so that
one knows how to deal with \eqref{3.6.1} in order to get
an {\em a priori\/} estimate of $D^{2}v$ on the set where
\eqref{3.6.1} holds through the values of $D^{2}v$ on the boundary
of this set, where \eqref{3.6.2} holds. However, one cannot justify
this line of arguments because one has to have a sufficiently
smooth $v$ from the very beginning.

Therefore, the argument in \cite{Kr12} uses finite-difference
approximations of \eqref{3.4.4}, but first one has to rewrite it
so that only pure second-order derivatives with respect to
a fixed (rather large but finite) family of vectors are
involved. One needs this because after replacing 
  second-order derivatives with  second-order differences
one wants to get a monotone finite-difference equation
(that is, an equation for which the maximum principle is valid).
We know how to do that, preserving
ellipticity, monotonicity in $u'_{0}$,
 and measurability or any kind of continuity
with respect to $x$,
 only if $H$ is boundedly inhomogeneous with respect to
$u,Du,D^{2}u$. In the situation of \cite{Kr13}
one can just replace $H$ with a different one,
which is boundedly inhomogeneous with respect to
$u,Du,D^{2}u$, and {\em without\/} changing  equation \eqref{3.4.4}.

After obtaining a finite-difference scheme one, basically,
repeats the above argument involving 
\eqref{3.6.2} and \eqref{3.6.1} dealing with the second-order
finite differences instead of $D^{2}v$. After
proving that the finite-difference versions of
\eqref{3.4.4} have solutions with the second-order
finite differences bounded independently of the mesh size,
we pass to the limit, as the mesh size goes to zero.
Once a sufficiently good solution of \eqref{3.4.4}
is secured, (a priori) $W^{2}_{p}(\Omega)$-estimates for that 
solution become available
owing to a different analytical line of arguments.
The estimates turned out to be independent of $K$
and sending $K\to\infty$ finishes the job.
 
In our present situation the auxiliary equation becomes 
\begin{equation}
                                                 \label{3.6.4}
\max(H(v( x),D v( x),D^{2}v( x), x),P( D^{2}v)-K)=0,
\end{equation}
because we can  hope to rewrite $H$ as a boundedly
inhomogeneous function only with respect to $u''$. This
turns out to be possible but the measure of inhomogeneity
becomes depending on $u'$, which 
seems  not to allow to preserve the {\em monotonicity\/}
in $u'_{0}$ of the corresponding finite-difference equations
and also led to other complications in estimates. 
In particular, in our heuristic argument, we cannot conclude
that $|D^{2}v|$ is bounded on the set where
$$
H(v( x),D v( x),D^{2}v( x), x)\geq P( D^{2}u)-K
$$
since there is no control of $|Du|$, and we cannot get it
by differentiating \eqref{3.6.4}. Furthermore,
while speaking about \eqref{3.4.4} above, 
 we said that we estimate $D^{2}v$, where
the opposite to \eqref{3.6.2} holds, through
the values of $D^{2}v$ on the boundary of this set,
to simplify the matter, we did not mention that
the boundary maybe part of $\partial\Omega$,
where there is no hope to get estimates for
$D^{2}v$ or for second-order differences.

In contrast with a parabolic version of Theorem
\ref{theorem 8.29.1}, presented in \cite{Kr13.1},
where losing monotonicity with respect to $u'_{0}$,
generally should not cause much problems (we only get our 
constants depending on the time length
of the cylinder in which the problem is  considered), in the
elliptic case it could just break the whole argument,
and an additional way to circumvent this
difficulty was needed.
 
The article is organized as follows. Sections
\ref{section 8.1.1} through \ref{section 8.29.1}
are devoted to the proof of Theorem \ref{theorem 8.29.1}.
In Section \ref{section 8.1.1} we present some existence results
about the solvability of nonlinear
 elliptic {\em finite-difference\/} schemes written
in terms of pure finite-differences.
In Section \ref{section 9_28.1} we cite
a particular case of Theorem 5.1 of \cite{Kr12},
yielding estimates for solutions of special nonlinear
finite-difference schemes with {\em constant\/}
coefficients. We apply them in Section \ref{section 9.21.1}
to obtain estimates for finite differences
for solution of ``cut-off" finite-difference equations (similar to
\eqref{9.23.02}). Section \ref{section 8.18.1} contains
existence results and estimates for ``cut-off'' {\em differential\/}
equations written in terms of pure second-order
derivatives. In Section \ref{section 8.25.1}, which is
central in the paper,
we consider general equations of ``cut-off'' type,
but with $H$ which is of ``restricted bounded inhomogeneity''
(in terms of \cite{Kr11})
and is Lipschitz with respect to $u'$. 
In Section \ref{section 8.29.1} we achieve the proof
of Theorem \ref{theorem 8.29.1} and in Section
\ref{section 3.14.1} we give the proof
of Theorems  \ref{theorem 2.20.1}
and \ref{theorem 2.20.2}. Finally, 
Section \ref{section 3.19.1} is an appendix
where, for the convenience of the reader,
 we present some well-known properties of Lipschitz
continuous functions.

\mysection{Solvability of elliptic finite-difference
equations}
                                          \label{section 8.1.1}

Fix a constant  
$m\in\{1,2,...\}$ and  
  a bounded domain $\Omega\subset\bR^{d}$, and 
for $h>0$ define
$$
\Omega^{h}=\{x:\rho_{\Omega}(x)>h\}
,\quad
\partial_{h}\Omega=\bar \Omega\setminus \Omega^{h},
$$
and let $\mu(\Omega)$ be the sup of $h$ such that $\Omega^{h}
\not=\emptyset$.
\begin{assumption}
                                   \label{assumption 7.31.1}   
We are given vectors $l_{i}\in \bar B_{1}\subset \bR^{d}$, $i=\pm1,...,
\pm m$,  
such that $l_{-i}=-l_{i}$ and $\Span\Lambda
=\bR^{d}$, where $\Lambda=\{l_{ i},|i|=1,...,m\}$.
\end{assumption}

For $h>0$ and $l\in\bR^{d}$ define $T_{h,l}\phi(x)=\phi(x+hl)$,
$$
\delta_{h,l}\phi(x)=\frac{1}{h}[T_{h,l}\phi(x)-\phi(x)].
\quad \Delta_{h,l}\phi(x)=\frac{1}{h^{2}}[\phi(x+hl)-2\phi(x)+
\phi(x-hl)].
$$

\begin{remark}
                                              \label{remark 8.4.1}
The reason to use $-l_{i}$ along with $l_{i}$ 
comes from the fact that in our calculations $-l_{i}$ would appear
anyway. 
\end{remark}

\begin{assumption}
                                \label{assumption 8.1.1}
(i\,) We are given a function $H(z,x)$, $ x\in\bar \Omega $,
$$
z=(z',z''),\quad z'=(z'_{0},z'_{\pm1},
...,z'_{\pm m})\in\bR^{2m+1},\quad z''=(z''_{\pm1},...,
z''_{\pm m})\in\bR^{2m},
$$
 which is Lipschitz continuous with respect to
$z$ for any $x$. At all point of its differentiability
with respect to $z$ introduce
$$
a_{ j}=a_{ j}(z,x)=D_{z''_{  j}}H(z,x),\quad
b_{ j}=b_{ j}(z,x)=D_{z'_{  j}}H(z,x),\quad j=\pm1,...,\pm m,
$$
$$
c =c(z,x)=-D_{z'_{0}}H(z,x),
$$
and at points of non-differentiability set $a_{ j}=\delta$,
$b_{ j}=0$, $c=1$. 

(ii\,) The above introduced functions
and $ H(0,x)$ are bounded,
 $|b_{j}|\leq N'$, $|j|=1,...,m$,  
where $N'$ is 
  a constant,  and
\begin{equation}
                                                   \label{8.1.1}
\delta^{-1}\geq a_{j}\geq \delta 
\end{equation}
for all indices and values of the arguments.

(iii\,) We have $c\geq0$.
\end{assumption}

For any function $v$ on $\bR^{d}$ and $h>0$ define
$$
H_{h}[v](x)=H(v,\delta_{h}v,\delta^{2}_{h} v,x)
$$
\begin{equation}
                                                      \label{8.8.2}
=H(v,\delta_{h,l_{\pm 1}}v,...,
\delta_{h,l_{\pm m}}v,\Delta_{h,l_{\pm 1}}v,...,
\Delta_{h,l_{\pm m}}v,x),
\end{equation}
where
\begin{equation}
                                                    \label{9.3.7}
\delta_{h}v=
(\delta_{h ,l_{\pm1}}u ,..., \delta_{h,l_{\pm m}}u ),
\quad \delta^{2}_{h} v=(\Delta_{h, l_{\pm1}}u ,..., \Delta_{h,l_{\pm m}}u ).
\end{equation}

Fix  
a bounded function $g$ on $\bar \Omega$ and consider
the equation
\begin{equation}
                                                 \label{8.3.1}
H_{h}[v]=0\quad\text{in}\quad \Omega^{h}
\end{equation}
with boundary condition
\begin{equation}
                                                 \label{8.3.2}
v=g\quad\text{on}\quad \partial_{h}\Omega.
\end{equation}

Of course, we will only consider $h$ such that $\Omega^{h}
\ne\emptyset$ ($h\in(0,\mu(\Omega)$). 
The following simple result
can be found in \cite{KT} or   in \cite{Kr11} 
or else in \cite{Kr12.1} if $c\geq1$. The proofs there are based on the method
of successive iteration, and the boundedness of $H(0,x)$
is used in order to start this procedure.
\begin{theorem}
                                         \label{theorem 8.3.1}
Let the above assumptions be satisfied. Then
for all sufficiently small $h>0$  the problem
\eqref{8.3.1}-\eqref{8.3.2} has a unique bounded solution
$v_{h}$.
\end{theorem}

To prove Theorem \ref{theorem 8.3.1} in the general
case we need the following, where
$$
D_{l_{j}}=l^{i}_{j }D_{i}.
$$

\begin{lemma}
                                           \label{lemma 8.3.1}
There exists a  function  $ \Psi_{0}
\in C^{\infty}(\bR^{d})$  
 such that, in $\Omega$,    $\Psi_{0}\geq1$ and
\begin{equation}
                                        \label{8.2.8}
a_{j} D^{2}_{l_{j}}\Psi_{0} +  
\delta^{-1}\sum_{j}|D_{l_{j}}\Psi_{0}|\leq-2,
\end{equation}
  and,
for  all sufficiently small $h>0$, in $\Omega^{h}$
\begin{equation}
                                        \label{8.3.7}
 a_{j}\Delta_{h,l_{j}}\Psi_{0} + \delta^{-1} 
\sum_{j}|\delta_{h,l_{j}}\Psi_{0} |\leq-1
\end{equation}
whenever $ \delta ^{-1}\geq a_{j}\geq \delta$.
\end{lemma}

The proof of this lemma is an elementary exercise. Indeed,
one takes 
$$
\Psi_{0}(x)=1+\cosh(\mu R)-\cosh(\mu|x|),
$$
where $R$ is twice the radius of the smallest ball 
centered at the origin containing
$\Omega$
and by straightforward computations using that
$\Span\{l_{j}\}=\bR^{d}$ and $a_{j}>\delta$ one gets
that (cf.~the proof of Lemma \ref{lemma 9.13.01})
$$
a_{j} D^{2}_{l_{j}}\Psi_{0} +  
\delta^{-1}\sum_{j}|D_{l_{j}}\Psi_{0}|\leq-\mu^{2}\cosh(\mu|x|)[\varepsilon-\mu^{-1}\delta^{-1}2m]
\leq-2 ,
$$
in $\Omega$ if $\mu$ is large enough, where $\varepsilon=\varepsilon(\Lambda)>0$.
This yields \eqref{8.2.8}. 

Then one observes that \eqref{8.2.8}   
is almost preserved if we replace $D_{l_{j}}$ and $D^{2}_{l_{j}}$
with $\Delta_{h,l_{j}}$ and $\delta_{h,l_{j}}$, respectively,
and choose $h>0$ sufficiently small. This is true
owing to the continuity of the derivatives of 
 $\Psi_{0}$
 in $\bR^{d}$.  

{\bf Proof of Theorem \ref{theorem 8.3.1}}.   
Take $\Psi_{0}$ from Lemma
 \ref{lemma 8.3.1}, in which we replace $\delta$
with a smaller one $\delta_{1}$ in order to have 
$|b_{j}|\leq \delta_{1}^{-1}$, and for $h>0$ set
$$
\bar H^{h}(z,x)=H( Z'_{0}(z,x),Z'_{\pm1}(z,x)...,Z'_{\pm m}(z,x),
 Z''_{\pm1}(z,x),...,Z''_{\pm m}(z,x),x),
$$
where, $Z'_{0}(z,x)=z'_{0}\Psi_{0}(x)$ and, for $|j|=1,...,m$,
$$
Z'_{j}(z,x)=z'_{j}T_{h,l_{j}}\Psi_{0}(x)+z'_{0}\delta_{h,l_{j}}\Psi_{0}(x),
$$
$$
Z''_{j}(z,x)=z''_{j}\Psi_{0}(x)+ z'_{j}\delta_{h,l_{j}}\Psi_{0}(x)+
z'_{-j}\delta_{h,l_{-j}}\Psi_{0}(x)+
z'_{0}\Delta_{h,l_{j}}\Psi_{0}(x)
$$
(no summation in $j$).
Note the fundamental property of $\bar H^{h}$
which follows by simple arithmetics: for any function $u$,  
\begin{equation}
                                                      \label{8.3.4}
\bar H^{h}_{h}[u]=H_{h}[u\Psi_{0}].
\end{equation}
In addition, if we denote by $(\bar a^{h} _{j},\bar b^{h}_{j},
\bar c^{h})(z,x)$
the functions corresponding to $\bar H^{h}$ as in Assumption 
\ref{assumption 8.1.1}, then obviously the boundedness
and nondegeneracy conditions in Assumption 
\ref{assumption 8.1.1} will be satisfied, and, due to Lemma
\ref{lemma 8.3.1},
$$
\bar c^{h} = c\Psi_{0}-b_{j}\delta_{h,l_{j}}\Psi_{0}
-a_{j}\Delta_{h,l_{j}}\Psi_{0}\geq1,
$$
where we dropped obvious values of the arguments for simplicity.

By what was said before the theorem, for sufficiently small $h>0$, the equation 
$\bar H^{h}_{h}[u]=0$ in $\Omega^{h}$ with boundary
condition $u=g/\Psi_{0}$ has a unique bounded solution.
It only remains to set $v=u\Psi_{0}$ and use 
\eqref{8.3.4}. The theorem is proved.  

\begin{lemma}[comparison principle]
                                           \label{lemma 8.4.1}
If $h>0$ is sufficiently small, then for any
  bounded functions 
$u,v$ on $\bar\Omega$
such that $u\leq v$ on $\partial_{h}\Omega$
and $H_{h}[u]\geq H_{h}[v]$ in $\Omega^{h}$
we have $u\leq v$ in $\bar\Omega$.
\end{lemma}

If $c(z,x)\geq\varepsilon>0$, the lemma is a particular case 
of Theorem 2.2 of \cite{Kr12.1}. 
In the general case it suffices to use the argument
in the end of the proof of Theorem \ref{theorem 8.3.1}.

In the future we will need an estimate of $v_{h}$.
\begin{theorem}
                                      \label{theorem 9.3.01}
Fix a constant $K_{0}\in[0,\infty)$ and
in addition to Assumptions \ref{assumption 7.31.1}
 and \ref{assumption 8.1.1} suppose that
 the number 
$$
\hat{ H}=\sup_{ x\in\Omega,z':z'_{0}=0 }(| H( z' ,0, x)|
-K_{0}| z' |)
$$
is finite. Then there is a constant $N$, depending only on  
 $K_{0}$, $\delta$, $\Lambda$,
and $\diam(\Omega)$,
  such that for  sufficiently small $h>0$
$$
\sup_{ \bar\Omega }|v_{h}|\leq N(\hat{ H}+\sup_{ \bar\Omega }|g |) .
$$

\end{theorem}

\begin{remark}
 Under Assumption \ref{assumption 8.1.1},
owing to the Lipschitz continuity of $H$ with respect to $z$
(uniform with respect to $x$ in light of
Assumption \ref{assumption 8.1.1} (ii)),
there always exists a $K_{0}$ such that
$\bar H<\infty$.
In the future, however, it will be important not to 
tighten up $K_{0}$ to the Lipschitz continuity.
\end{remark}

The proof of Theorem \ref{theorem 9.3.01} is based on the following.

\begin{lemma}
                                            \label{lemma 9.3.10}
Under the assumptions of Theorem \ref{theorem 9.3.01}
let $\phi$ be a function on $\Omega$. Then
for any $h\in(0,\mu(\Omega))$  \(that is such that $\Omega^{h}\ne\emptyset$\),
 there exist bounded functions $a_{j},b_{j},c,f$,
$|j|=1,...,m$, in $\Omega^{h}$
\(generally different from the ones
in Assumption \ref{assumption 8.1.1}\)
 such that in $\Omega^{h}$ 
$$
\delta \leq a_{j}\leq  \delta^{-1},\quad |b_{j}|\leq K_{0},
\quad c\geq0,\quad |f|\leq\hat H,
$$
\begin{equation}
                                                       \label{9.3.5}
  H(\phi,\delta_{h}\phi ,
\delta_{h}^{2}\phi, x) 
=a_{k}\Delta_{h,l_{k}}\phi +
b_{k}\delta_{h,l_{k}}\phi-c\phi+f.
\end{equation}

\end{lemma}

Proof.
Notice that, due to Assumption \ref{assumption 8.1.1}
(cf.~Section \ref{section 3.19.1}),
$$
  H(\phi,\delta_{h}\phi ,
\delta_{h}^{2}\phi, x)  =  H(\phi,\delta_{h}\phi,
\delta_{h}^{2}\phi, x)  -  H (\phi,\delta_{h}\phi,
0, x)
$$
$$
 +H (\phi,\delta_{h}\phi,
0, x) 
  = a_{k}\Delta_{h,l_{k}}\phi 
+H (\phi,\delta_{h}\phi,
0, x) ,
$$
where $a_{k}$ are some functions satisfying $
\delta \leq a_{k}\leq \delta^{-1}$,
and 
$$
  H (\phi,\delta_{h}\phi,
0, x) =-c\phi+f
$$
with bounded $c\geq0$ and 
$f:=H(0,\delta_{h}\phi,0, x )$ satisfying
$$
|f|\leq 
\hat{ H} +K_{0}  \sum_{|k|=1}^{m}
|\delta_{h,l_{k}}\phi|   .
$$

This property of $f$ implies that there exist
functions
$b_{k}$, $|k|=1,...,m$,  
with values in $[-K_{0},K_{0}]$  and $\theta$
with values in $[-1,1]$
 such that
$$
f( \delta_{h}\phi, x)=  
b_{k}\delta_{h,l_{k}}\phi+\theta\hat{  H} .
$$
Upon combining all the above we come to \eqref{9.3.5}.
The lemma is proved.  

{\bf Proof of Theorem \ref{theorem 9.3.01}}.
In Lemma \ref{lemma 8.3.1} reduce $\delta$, if necessary,
in such a way that the new one, say $\delta_{1}$,
  satisfy $K_{0}\leq\delta^{-1}_{1}$,
and then take $\Psi_{0}$ from that lemma
and set
$$
\phi=\Psi_{0}(\hat{ H}+\sup_{ \bar\Omega }|g |).
$$
Then by using Lemmas \ref{lemma 9.3.10} and  \ref{lemma 8.3.1}
we see that for   $h\in(0,h_{0}]$ on $\Omega^{h}$
$$
 H(\phi,\delta_{h}\phi ,
\delta_{h}^{2}\phi, x)=a_{k}\Delta_{h,l_{k}}\phi
+b_{k}\delta_{h,l_{k}}\phi
-c\phi+\theta\hat{  H}\leq-(\bar{ H}+\sup_{ \bar\Omega }|g |)
+\theta\hat{  H}\leq0.
$$
Furthermore, $ v^{h} = g \leq\phi$ on $\partial_{h}\Omega$.
 It follows by Lemma \ref{lemma 8.4.1}
that $v^{h}\leq\phi$ on $\bar\Omega$. By replacing $\phi$
with $-\phi$ we get that $v^{h}\geq-\phi$ on $\bar\Omega$.
The theorem is proved.

\mysection 
{Some estimates for finite-difference  
equations with constant coefficients}
                                            \label{section 9_28.1}
Take an $h\in(0,1]$, let
$m\geq1$ be an integer and let
$l_{\pm1},...,l_{\pm m} $ be some fixed vectors in
$\bR^{d}$ such that  $ l_{-k}=-l_{k}$. 
Denote
$$
\Lambda=\{l_{k}:k=\pm1,...,\pm m\},
$$  
$$ \Lambda_{1}= \Lambda,\quad
\Lambda_{n+1}=
 \Lambda_{n}+ \Lambda ,\quad n\geq1,
\quad \Lambda_{\infty}=\bigcup_{n}\Lambda_{n},
\quad \Lambda^{h}_{\infty}=h\Lambda_{\infty}
\,. 
$$
 
Let
$Q^{o}$ be  a nonempty finite subset
of $  \Lambda_{\infty}^{h}$.
Introduce  
\begin{equation}
                                                 \label{6.11.1}
Q=Q^{o}\cup\{ x+h\Lambda : x \in Q^{o}\}.
\end{equation}
  
Let   $A $ be a
closed bounded set of points
$$
a=(a_{\pm1}, ...,a_{\pm m})\in\bR^{2m}.
 $$
\begin{assumption}
                                   \label{assumption 9.18.1}
 There is a constant  $\delta\in(0,1]$  such that
for any $a\in A$ and all $k$ we have  $ a_{k}=a_{-k}$ and
$\delta \leq a_{k}\leq \delta^{-1}$.

\end{assumption}
 Also let $f(a )$
be a real-valued continuous function defined on $  A$ and for  
$$
 z'' =(z''_{\pm1},...,z''_{\pm m} )\in\bR^{2m},
$$
introduce
$$
 P(z'')=\max_{a\in A}\big(\sum_{|k|=1}^{m}
a_{k}z''_{k}  +f(a)\big).
$$

For any function $u$ on $\bR^{d }$ define
\begin{equation}
                                                               \label{9.27.1}
  P_{h}[u]( x)= P(  \delta^{2}_{h}u( x)  ),
\end{equation}
where $\delta^{2}_{h}$ is introduced in \eqref{9.3.7}.

In connection with this notation a
natural question arises as to why use $l_{k}$ along with
$l_{-k}=-l_{k}$ since $\Delta_{h,l_k}=\Delta_{h,l_{-k}}$ and 
$$
a_{k}\Delta_{h,l_k}=2\sum_{k\geq1}a_{ k}\Delta_{h,l_k} 
$$ 
owing to the
assumption that $a_{k}=a_{-k}$. This is done for the sake of convenience
of computations. For instance,
$$
\Delta_{h,l_k}(uv)=u\Delta_{h,l_k}v+
v\Delta_{h,l_k}u+(\delta_{h,l_k}u)(\delta_{h,l_k}v)
+(\delta_{h,-l_k}u)(\delta_{h,-l_k}v)
$$ (no summation in $k$). At
the same time
$$ a_{k}\Delta_{h,l_k}(uv)=ua_{k}\Delta_{h,l_k}v+va_{k}
\Delta_{h,l_k}u
+2a_{k}(\delta_{h,l_k}u)(\delta_{h,l_k}v),
$$
as if we were dealing with usual
partial derivatives. Another, even more compelling, reason
is mentioned in Remark \ref{remark 8.4.1}.

Next,
take a function $\eta\in C^{\infty}(\bR^{d})$ with bounded
derivatives, such that $|\eta|\leq1$ and set
$\zeta=\eta^{2}$, $$   |\eta'(x)|_{h} =\sup_{
k}|\delta_{h,l_k}\eta  (x)|,\quad   |\eta''(x)|_{h}
=\sup_{ k}|
\Delta_{h,l_k}\eta  (x)|, $$ $$
 \|\eta'\|_{h}=\sup_{ \Lambda_{\infty}^{ h }}|\eta'
|_{h},\quad
 \|\eta''\|_{h}=\sup_{
\Lambda_{\infty}^{ h}}|\eta''
|_{h}. 
$$

Finally, let $u$ be a function on $\bR^{d }$ which   satisfies
\begin{equation}
                                                \label{9.19.9}
 P_{h}[u]=0 \quad\text{in}\,\,Q^{o}
\end{equation}
 and
\begin{equation}
                                        \label{eq4.48}
 P_{h}[u]\leq 0\quad
\text{on}\,\, Q \setminus Q^{o} . 
\end{equation}

Here is a one-sided interior estimate of pure second-order
differences of~$u$.

\begin{theorem}
                                      \label{theorem 9.5.1}
There  exists a constant  $N=N(m,\delta)\geq1$
such that for any $r=\pm1,...,\pm m$ \(recall
that  $a^{\pm}=(1/2)(|a|\pm a)$\) in $Q$ we have
\begin{equation}
                                               \label{9.5.7}
\zeta^{2}[ ( \Delta_{h,l_r}u)^{-}]^{2}\leq \sup_{
 Q\setminus Q^{o} }\zeta^{2}[
( \Delta_{h,l_r}u)^{-}]^{2} + N(\|\eta''\|_{h}+\|\eta'\|_{h}^{2}) \bar{W}_{r},
\end{equation}
where
\begin{equation}
                                                      \label{9.5.1}
\bar{W}_{r}=\sup_{ Q  }(|\delta_{h,l_r}u|^{2}+
|\delta_{h,l_{-r}}u|^{2}).
\end{equation}

\end{theorem}

This theorem is a particular case of Theorem 5.1 of \cite{Kr12}.
 
\mysection{Estimates for
finite-difference   equations
of cut-off type with variable coefficients}
                                            \label{section 9.21.1}

Fix some constants  $K_{0},K\in[0,\infty)$.
We use the notation and assumptions introduced   
in the beginning of Section
\ref{section 8.1.1}, however, we append Assumptions  
\ref{assumption 7.31.1} and
\ref{assumption 8.1.1} with the following.  

\begin{assumption}
                                   \label{assumption 8.10.1}
The coordinates of $l_{k},|k|=1,...,m$, are rational numbers.
\end{assumption}

\begin{assumption}
                                   \label{assumption 9.23.01}
 The number 
$$
\bar{ H}=\sup_{ x,z' }(| H( z' ,0, x)|
-K_{0}| z' |)
$$  
is finite.
\end{assumption}

 We also impose the following.   

\begin{assumption}
                                   \label{assumption 8.8.1}

The domain $\Omega$ is bounded and satisfies the exterior
ball condition.
 We are given a function $g\in C^{1,1}(\bR^{d})$.
\end{assumption}

For $z''\in\bR^{2m}$ introduce
\begin{equation}
                                               \label{6.1.1}
 P(z'')  =\max_{\substack{ \delta /2\leq a_{k}\leq
2 \delta ^{-1} \\|k|=1,...,m} }
 \sum_{|k|=1}^{m} a_{k}
z''_{k} 
\end{equation}

Observe that
\begin{equation}
                                                       \label{9.27.2}
 H(z, x)\leq  P(z'')-(\delta/2)\sum_{|k|=1}^{m}|z''_{k}|
+K_{0}|z'|+\bar{ H}.
\end{equation}
Indeed,   it follows from Assumption \ref{assumption 8.1.1}  that 
(cf.~Section \ref{section 3.19.1})
$$
 H(z,x)- H(z',0, x)=\sum_{|k|=1}^{m}a_{k}z''_{k},
$$
where $\delta\leq a_{k}\leq\delta^{-1}$ for all $k$. Hence,
$$
 H(z,x)- H(z',0, x)+(\delta/2)\sum_{|k|=1}^{m}|z''_{k}|
=\sum_{|k|=1}^{m}b_{k}z''_{k},
$$
where $b_{k}:=a_{k}+(\delta/2)\sign z''_{k}$ lie
in $(\delta/2,\delta^{-1}+\delta/2)\subset(\delta/2,2/\delta)$.

We need one more function
$$
H_{K}(z,x)=\max(H(z,x),P(z'')-K).
$$
Recall that $\delta_{h}^{2}v$ and $\delta_{h}v$ are defined in 
\eqref{9.3.7}   and for any function $v$ on $\Omega$
set
$$
H_{K,h}[v]( x)
=H_{K}(v( x),\delta_{h}v( x),\delta_{h}^{2}v( x), x).
$$

We will concentrate on sufficiently small $h$
such that $\Omega^{  h}\ne\emptyset$ 
and will
consider the equation
\begin{equation}
                                              \label{2.25.3}
 H_{K,h}[v]=0\quad \text{in}\quad
 \Omega^{ h}   
\end{equation}
 with  boundary condition
\begin{equation}
                                              \label{2.25.4}
v=g\quad\text{on}\quad  \partial_{h}\Omega . 
\end{equation}

By Theorem \ref{theorem 8.3.1}  for any sufficiently small $h>0$ there
exists a unique bounded solution $v=v_{ h}$ of
\eqref{2.25.3}--\eqref{2.25.4}.
  By the way, we
do not include $K$ in the notation $v_{ h}$ since $K$ is a fixed
number. The solution of the PDE version 
of \eqref{2.25.3}--\eqref{2.25.4} will be obtained
as the limit of a subsequence of $v_{h}$ as $h\downarrow0$.
Therefore, we need to have appropriate bounds on
 $v_{h}$   and the first- and second-order differences in $x$
of $v_{h}$.

Observe that owing to Assumption \ref{assumption 9.23.01},  
for $z'_{0}=0$,
$$
- K_{0} |z'|-\bar H\leq H (0,z'_{\pm1},...,z'_{\pm m},0,x)
\leq H_{K}(0,z'_{\pm1},...,z'_{\pm m},0,x)
$$
$$
=\max(H (0,z'_{\pm1},...,z'_{\pm m},0,x),-K)\leq K_{0} |z'|+\bar H,
$$
so that $\hat H_{K}\leq\bar H$, where $\hat H_{K}$ is taken
from Theorem \ref{theorem 9.3.01} with $H_{K}$ in place of $H$.
Therefore, the following is a direct consequence of 
Theorem~\ref{theorem 9.3.01}.  
   
 \begin{lemma}
                                    \label{lemma 10.6.1}  
There is a constant $N$, depending only on $\Lambda$,   $K_{0}$, $\delta$,
and $\diam(\Omega)$,
 such that for sufficiently small $h>0$
\begin{equation}
                                                          \label{8.17.1}
\sup_{ \bar\Omega }|v_{h}|\leq N(\bar{ H}+\|g\|_{C(\Omega )}).
\end{equation}

\end{lemma}

Lemma \ref{lemma 9.3.10}   also
  yields the following
useful  result.

\begin{lemma}
                                            \label{lemma 8.8.2}
Let $\phi$ be a function on $\Omega$. Then
for any $h\in(0,\mu(\Omega))$ there exist bounded functions $a_{j},b_{j},c,f$,
$|j|=1,...,m$, on $\Omega^{h}$ such that on $\Omega^{h}$ 
\begin{equation}
                                                            \label{8.8.60}
\delta/2\leq a_{j}\leq 2\delta^{-1},\quad |b_{j}|\leq K_{0},
\quad c\geq 0,\quad |f|\leq\bar H,
\end{equation}
$$
  H_{K,h}[\phi]   
=a_{k}\Delta_{h,l_{k}}\phi +
b_{k}\delta_{h,l_{k}}\phi-c\phi+f.
$$

\end{lemma}

Below by   $N$ with occasional indices we denote various 
(finite) constants depending only on $\diam(\Omega)$,
$\rho(\Omega)$, $\mu(\Omega)$,  
$\Lambda=\{l_{i}\}$, $K_{0}$, 
$d$, and $\delta$, unless explicitly stated otherwise.

We need a barrier function. Recall that   $\rho(\Omega)\leq1$.
 \begin{lemma}
                                         \label{lemma 9.13.01}
There exists a constant $\alpha>0$, depending only on
$\Lambda$,  $\delta$, $\rho (\Omega)$, and $K_{0}$, such that for   
$\psi(x)=|x|^{-\alpha}-\rho^{-\alpha}(\Omega)$ we have
\begin{equation}
                                              \label{9.29.1}
  a_{k}D^{2}_{ l_{k} }\psi(x)
+ b_{k}D _{ l_{k} }\psi(x)-c\psi(x)\geq 1,
\end{equation}
whenever $2\delta^{-1}\geq a_{k}\geq\delta/2$, $|b_{k}|\leq K_{0}$, 
$|k|=1,...,m$, $|c|\leq K_{0}$,
and $0<|x|\leq 3$.
\end{lemma}

Proof. Observe that
$$
  a_{k}D^{2}_{ l_{k} }\psi =
a^{ij}D_{ij}\psi,
$$
where
$$
a^{ij}= a_{k}l_{k}^{i}l_{k}^{j}.
$$
For any $\lambda\in\bR^{d}$ it holds that
$$
a^{ij}\lambda_{i}\lambda_{j}=
 a_{k}(l_{k},\lambda)^{2}
\geq\delta
\sum_{k=1}^{m} (l_{k},\lambda)^{2}\geq\delta_{1}|\lambda|^{2},
$$
$$
a^{ij}\lambda_{i}\lambda_{j}\leq2m\delta^{-1} |\lambda|^{2}
=:N_{1}|\lambda|^{2},
$$
where $\delta_{1}>0$ is a constant whose existence
is guaranteed by the assumption that
$\Span\{l_{j}\}=\bR^{d}$.

For    $\alpha$ such that
$(\alpha+2)\delta_{1}\geq N_{1}d+\alpha\delta_{1}/2$ we obtain that
$$
a^{ij}D_{ij}\psi=\alpha|x|^{-\alpha-2}
((\alpha+2)|x|^{-2}a^{ij}x_{i}x_{j}-\tr (a^{ij}))
\geq (\delta_{1}/2)\alpha^{2}|x|^{-\alpha-2}.
$$
One easily finishes the proof after noticing that
$$
\big| b_{k}D_{l_k}\psi\big|=
\alpha|x|^{-\alpha-2}|b_{k}(l_{k},x)|,\quad|c\psi|\leq K_{0}
(|x|^{-\alpha}+\rho^{-\alpha}(\Omega)).
$$
The lemma is proved.

 \begin{lemma}
                                         \label{lemma 9.13.1}
There   is a
constant  $N$
 such that, for all sufficiently small $h>0$,
\begin{equation}
                                              \label{9.13.1}
|v_{ h}-g|\leq N(\bar{ H}+ 
\|g\|_{C^{1,1}( \Omega )  })(\rho_{\Omega}\wedge 1)
\end{equation}
on $ \bar\Omega $.
\end{lemma}

Proof. 
Observe that \eqref{2.25.4} implies that
\eqref{9.13.1} holds in $\partial_{h}\Omega  $. 
Therefore, we may concentrate on proving
\eqref{9.13.1} in $ \Omega^{  h}$.

Take  an $x_{0}\in \Omega^{  h}$
(so that $  h<\rho_{\Omega}(x_{0})$). If $\rho_{\Omega}(x_{0})\geq1$,
then \eqref{9.13.1} holds at $ x_{0} $ in light of
Lemma \ref{lemma 10.6.1}. Hence we may assume that
$\rho_{\Omega}(x_{0})\leq1$. 
Let $y$ be a closest point to $x_{0}$ on
$\partial \Omega$. Then the closures of the following two balls
 has only one common point $y$: one is
the ball centered
at $x_{0}$ with radius $\rho_{\Omega}(x_{0})$
and the other is the  ball of radius $\rho(\Omega)$,
which lies outside $\Omega$ and is centered, say at $y_{0}$.
 Then all three points
$x_{0}$, $y$, and $y_{0}$ lie on the same line.
Therefore,  without loss of generality, 
on the account of moving the origin,  we may
assume that $y_{0}=0$,
the straight line passing through the origin and $x_{0}$
crosses $\partial \Omega$ at $y$, and $|y|= \rho(\Omega)  $ ($\leq1$
by assumption)
 
Then set
$$
\phi=g-2N_{1}(\bar H+\|g\|_{C^{1,1}(\Omega)}) \psi ,
$$
where  $\psi$ is taken from Lemma \ref{lemma 9.13.01}
and the constant $N_{1}$ is such that in $\Omega^{h}$
$$
a_{k}\Delta_{h,l_{k}}g +
b_{k}\delta_{h,l_{k}}g-cg+f\leq  
N_{1}(\bar H+\|g\|_{C^{1,1}(\Omega)})  
$$
as long as
conditions \eqref{8.8.60} are satisfied.

We increase $N_{1}$ if necessary in order to have
$v_{h}\leq \phi$ in $\Omega\setminus B_{2}$,
 the latter being possible owing to
Lemma \ref{lemma 10.6.1} and the fact that 
in $\Omega\setminus B_{2}$
we have 
$$
 -\psi\geq 
\rho^{-\alpha}(\Omega)-2^{-\alpha} 
\geq 1-2^{-\alpha}.
 $$
 We also concentrate on $h>0$ sufficiently small such
 that the finite-difference approximations
of $D^{2}_{ l_{k} }\psi$ and $D_{ l_{k} }\psi$ are sufficiently
close to these quantities in $B_{3}\setminus B_{\rho(\Omega)}$,
so that owing to Lemma \ref{lemma 9.13.01}
$$
 a_{k}\Delta_{h,l_{k}}\psi 
+ 
b_{k}\delta_{h,l_{k}}\psi-c\psi\geq 1/2
$$
in $B_{3}\setminus B_{\rho(\Omega)}$.

Then for $\Omega_{1}=\Omega\cap B_{3}$
(where $ \psi<0$), by using Lemma 
\ref{lemma 8.8.2}, we find that in $\Omega_{1}^{h}$
$$
\max\big( H(\phi,\delta_{h}\phi ,
\delta_{h}^{2}\phi, x),P_{h}[\phi]-K\big)\leq0.
$$
Also by the choice of $N_{1}$,
$v_{h}\leq \phi$ on $\partial_{h}\Omega_{1}$  (if $h<1$).
Hence by the comparison principle, $v_{h}
\leq \phi$ in $\Omega_{1}=\Omega\cap B_{3}$.
By the choice of $N_{1}$, this inequality,
actually, holds in $\Omega$, 
which yields the desired estimate of
$v_{h}-g$ from above. Similarly one obtains it from below as well. 
The lemma is proved. 

The proof of the lemma allows us 
to get control on the boundary behavior of $v_{h}$
uniformly with respect to $h$ in a relatively easy way. The way to treat 
finite-differences of $v_{h}$ is much more involved.

For sufficiently small $h>0$ (such that $v_{h}$
is well defined) introduce
$$
M_{h}( x)=\sum_{|k|=1}^{m}|\delta_{h,l_{k}}
v_{h}( x)|,\quad \bar M_{h}=\sup_{\Omega^{h}}M_{h},
$$ 
\begin{equation}
                                             \label{9.25.1}
G:=\{ x \in\Omega^{  h} :  (\delta/2)  \sum_{k=1}^{m}
|\Delta_{h,l_{k}}v_{h}( x)|
>\bar{ H}+K+K_{0}
\big(|v_{h}( x)|+M_{h}( x)  
\big)\}.
\end{equation}
and observe a fundamental fact that thanks to 
\eqref{9.27.2}  (cf.~\eqref{3.6.1})
\begin{equation}
                                              \label{9.16.7}
 P_{h}[v_{h}]-K=0 
\quad \text{in}\quad G.
\end{equation}

Also observe that on $\Omega^{h}\setminus G$ we have
\begin{equation}
                                                   \label{9.22.7}
(\delta/2)  \sum_{k=1}^{m}
|\Delta_{h,l_{k}}v_{h} |
\leq\bar{ H}+K+K_{0}
\big(|v_{h} |+M_{h}  ) .
\end{equation}

\begin{lemma}
                                          \label{lemma 9.16.1}
There  is a
constant  $N$  such that, for
all sufficiently small $h>0$ and $r=1,...,m$,  in
$ \Omega^{  h}$ we have
\begin{equation}
                                              \label{9.16.6}
  (\rho_{\Omega}-2h) |\Delta_{h,l_{r}}  v_{ h}|\leq N(
\bar M_{h}+
\bar{ H}+K
+\|g\|_{C^{1,1}( \Omega )}) .
\end{equation}
\end{lemma}

Proof.
Having in mind 
translations, we see that it suffices to prove \eqref{9.16.6} in
$ \Omega^{2h} 
\cap
\Lambda_{\infty}^{h} $.  Then  
  fix $r$ and define
$$
Q^{o}:=\Omega^{ 2 h} 
\cap
\Lambda_{\infty}^{h}\cap G.
$$
In light of Assumption \ref{assumption 8.10.1} (used for the first
time), the set $Q^{o}$ is finite, since the number of
points in $\Lambda^{h}_{\infty}$ lying in any ball is finite,
because for an appropriate integer $I$ we have $I\Lambda^{h}_{\infty}
\subset h\bZ^{d}$.

 For $Q$ from \eqref{6.11.1},  obviously,   
$Q\subset \Omega^{  h}$. Next,
if   $x\in \Omega^{2h}\cap\Lambda^{h}_{\infty}$ is such
that $x \not\in Q^{o}$, then  \eqref{9.22.7} is valid,
 in which case \eqref{9.16.6} holds.

Thus, we need only prove \eqref{9.16.6} on $Q^{o}$ assuming, of course,
that $Q^{o} \ne\emptyset$. We know that 
  \eqref{9.16.7}   holds and the left-hand side of
\eqref{9.16.7} is nonpositive  in $Q\setminus Q^{o}$
(as everywhere else in $\Omega^{h}$).  

To proceed further we use a simple fact that there exists
a constant $N\in(0,\infty)$ depending only on $d$ such that
for any
$\mu \in(0,
\mu(\Omega)/2)$
  there exists an $\eta_{\mu}
\in C^{\infty}_{0}(\Omega)$ satisfying
$$
\eta_{\mu}=1\quad\text{on}\quad \Omega^{2 \mu},\quad
\eta_{\mu}=0\quad\text{outside}\quad \Omega^{  \mu},
$$
$$
|\eta_{\mu}|\leq1,\quad |D\eta_{\mu}|\leq
N/\mu,\quad|D^{2}\eta_{\mu}|\leq N/\mu^{2}.
$$
By  Theorem \ref{theorem 9.5.1} in $Q^{o}  $
$$
 [ ( \Delta_{h,l_r}v_{h})^{-}]^{2}\leq \sup_{
 Q\setminus Q^{o} }\eta_{\mu}[
( \Delta_{h,l_r}v_{h})^{-}]^{2} + N
\mu^{-2}\bar M_{h}^{2}.
$$
While estimating the last supremum
we will only concentrate on (sufficiently small $h$ and)  
$$
 \mu\in[ 2h,\mu(\Omega)/2),
$$
 when
$\eta_{\mu}=0$ outside $\Omega^{2  h}$. In that case, for any
$y\in
  Q\setminus Q^{o}$, either $y\notin\Omega^{2  h}$ implying that
$$
\eta_{\mu} [ ( \Delta_{h,l_{r}}v_{h})^{-}]^{2}( y)=0,
$$
or  $y \in\Omega^{2  h}\cap \Lambda_{\infty}^{h}$ but 
\eqref{9.22.7} holds at $y$. 

It follows that, as long as   $x\in Q^{o}
  $  and  $\mu\in
[2h,\mu(\Omega)/2)$, we have   
\begin{equation}
                                                   \label{9.22.07}
\mu  ( \Delta_{h,l_{r}}v_{h})^{-}( x)\leq
 N
(\bar{ H}+K
+\|g\|_{C^{1,1}( \Omega )}+\bar M_{h}).
\end{equation}
Take $$
\mu=(\mu(\Omega)/2)\wedge \rho_{\Omega}(x) ,
$$
 which is bigger than $2h$
on $Q^{o}$. Also $\mu=\rho_{\Omega}(x)$
if $\rho_{\Omega}(x)\leq \mu(\Omega)/2 $ and
$$
\mu= \mu(\Omega)/2 \geq \rho_{\Omega}(x) \mu(\Omega)/(2\,\diam(\Omega) )
$$
if $\rho_{\Omega}(x)\geq \mu(\Omega)/2 $. This and \eqref{9.22.07}
yield that on $Q^{o}$
\begin{equation}
                                                      \label{10.9.1} 
\rho_{\Omega}( \Delta_{h,l_{r}}v_{h})^{-}( x)
\leq N (\bar{ H}+K
+\|g\|_{C^{1,1}( \Omega )}+\bar M_{h}).
\end{equation}
Obviously, one can replace $\rho_{\Omega}$ in \eqref{10.9.1}   with $ \rho_{\Omega}-3 h $  and as a result of all the above arguments we see that 
\begin{equation}
                                                   \label{9.22.8}
(\rho_{\Omega}-3 h)( \Delta_{h,l_{r}}v_{h})^{-}\leq N(
\bar{ H}+K
+\|g\|_{C^{1,1}( \Omega )}+\bar M_{h}) 
\end{equation}
 holds
in $Q^{o}$
  for any $r$ whenever $h $ is small enough.  

Finally, since $ P_{h}[v_{h}]\leq K$ in $  \Omega^{  h}$, we have that
$$
4 \delta ^{-1}\sum_{r=1}^{m}(\Delta_{r}v_{h})^{ +} 
\leq  \delta 
\sum_{r=1}^{m}(\Delta_{r}v_{h})^{ -}+K,
$$
which after being multiplied by $(\rho_{\Omega}-3 h)^{+}$ along
with \eqref{9.22.8}  
 leads to \eqref{9.16.6}
on $Q^{o}$. Thus, as
  is explained at the beginning of the proof,
  the lemma is proved.                                

Now we  exclude $\bar M_{h}$
from \eqref{9.16.6}.
 
\begin{lemma}
                                       \label{lemma 9.22.3}
 
There is a constant $N$  such that for all sufficiently
small $h>0$ the estimates
\begin{equation}
                                                  \label{9.18.9}
 |v_{h}|, |\delta_{h,l_{k}}v_{h}|,
 ( \rho_{\Omega} -2h) |\Delta_{h,l_{k}}v_{h}|\leq N
(\bar{ H}+K
+\|g\|_{C^{1,1}( \Omega )})
\end{equation}
hold  in $ \Omega^{  h}$ for all $k$.
\end{lemma}
 
Proof.  Owing to Lemmas \ref{lemma 9.13.1}
and \ref{lemma 9.16.1},   \eqref{9.18.9}
would follow if we can prove that
\begin{equation}
                                                   \label{9.22.08}
 |\delta_{h,l_{k}}v_{h}|\leq N
(\bar{ H}+K
+\|g\|_{C^{1,1}( \Omega )})
\end{equation}
   in $ \Omega^{  h}$ for all $k$. Actually, estimate
\eqref{9.22.08} is proved in Lemma 4.7 of \cite{Kr13.1}
for the parabolic case by using interpolation
inequalities, but the proof
is valid word for word for the elliptic case as well.
The lemma is proved.

\begin{remark}
                                    \label{remark 11.6.1}
Much of what is done above goes through
for domains satisfying the exterior cone
condition instead of the exterior ball condition.
However, then it would be impossible to get
the global discrete gradient estimate and to get
  estimates of $\Delta_{h,l_{k}}v_{h}$ in a closed form.
\end{remark}

\mysection{A particular case of elliptic
equations in pure derivatives}
                                      \label{section 8.18.1}

Fix some constants 
$m\in\{1,2,...\}$  and $K_{0},K, N'
 \in(0,\infty)$.
 
\begin{assumption}
                                  \label{assumption 8.18.1}
We are given vectors $l_{i}\in \bar B_{1}\subset \bR^{d}$, $i=\pm1,...,
\pm m$,
such that $l_{-i}=-l_{i}$, the coordinates
of $l_{i}$'s are rational numbers, 
and 
$$
e_{i}, \,(1/2)(e_{i}\pm e_{j})\in \Lambda:=\{l_{k}:k=\pm1,...,\pm m\}
$$
for $i,j=1,..,d$, $i\ne j$, where $e_{i}$'s form the standard
orthonormal basis in $\bR^{d}$.
\end{assumption}

\begin{assumption}
                                \label{assumption 8.18.3}

We are given a domain $\Omega
\subset\bR^{d}$ which is bounded and satisfies the exterior
ball condition. We are also
 given a function $g\in C^{1,1}(\bR^{d})$.
\end{assumption}

The following assumption has some parts which are
 close to Assumption
\ref{assumption 8.1.1}.

\begin{assumption}
                                \label{assumption 8.18.2}
(i) We are given a continuous 
function $H(z,x)$, $ x\in\bR^{d}$,
$$
z=(z',z''),\quad z'=(z'_{0}, z'_{\pm1},
...,z'_{\pm m})\in\bR^{2m+1},\quad z''=(z''_{\pm1},...,
z''_{\pm m})\in\bR^{2m},
$$
 which is Lipschitz continuous with respect to
$z\in \bR^{4m+1}$ with Lipschitz constant $N'$ for any $x\in\bR^{d}$. 
At all point of its differentiability
with respect to $z $ introduce
$$
a_{ j}=a_{ j}(z,x)=D_{z''_{  j}}H,\quad c=c(z,x)=-D_{z'_{0}}H 
$$
and at points of non-differentiability set $a_{ j}=\delta$
and $c=1$.

(ii) The above introduced functions $a_{j}$ and $c$
satisfy
\begin{equation}
                                             \label{8.1.01}
\delta^{-1}\geq a_{j}\geq \delta ,\quad c\geq0
\end{equation}
for all indices and values of the arguments.

(iii)  The number 
$$
\bar{ H}=\sup_{ x,z' }(| H( z' ,0, x)|
-K_{0}| z' |)
$$ is finite.

(iv) The function $H(z,x)$ is Lipschitz continuous
with respect to $x\in\bR^{d}$ for any $z\in\bR^{4m+1 } $ and at any
point $x$ of its differentiability with respect to~$x$
\begin{equation}
                                    \label{9.31.002}
|D_{x}H (z,x)|\leq N' (1+|z |).
\end{equation}

(v) There is a constant $\rho_{0} >0$ and a function $H(z )$ such that
\begin{equation}
                                            \label{8.17.02}
H(z,x)=H(z )
\end{equation}
for all $z$ if $\rho_{\Omega}(x)\leq\rho_{0} $
 \(recall that  
 $
\rho_{\Omega}=\rho_{\Omega}(x)=\dist(x,\bR^{d}\setminus\Omega)
 $\).

\end{assumption}

For sufficiently smooth functions $u=u(x)$
introduce
$$
H[u](x)=H(u(x), D_{l_{j}}u(x),D^{2}_{l_{k}}u(x),x)
$$
$$
=H(u(x), D_{l_{\pm1}}u(x),
...,D_{l_{\pm m}}u(x),
D^{2}_{l_{\pm1}}u(x),
...,D^{2}_{l_{\pm m}}u(x),x).
$$
Similarly we introduce $P[u]$, where
$P$ is taken from \eqref{6.1.1}.

 We will be
interested in finding
  a solution $v\in
C(\bar{\Omega} )\cap C^{1,1}_{ \loc }(\Omega)$ 
of equation \eqref{9.23.02} in $\Omega$ with boundary data $g$.
We will also be interested in obtaining the estimates
\begin{equation}
                                             \label{10.15.30}
 |v|\leq N(\bar H+\|g\|_{C(\partial\Omega)})
\quad\text{in} 
\quad \Omega,
\end{equation}
\begin{equation}
                                             \label{9.22.3}
|v|,|D_{i}v|, \rho_{\Omega} |D_{ij} v | \leq N(\bar{H}+K
+\|g\|_{ C^{1,1} (\bR^{d})})\quad\text{in} 
\quad \Omega  \quad \text{(a.e.)},
\end{equation}
where $i,j=1,...,d$,
 $D_{ij}v$ are generalized derivatives of $v$, and
the constants $\bar H$  and $N$
are chosen appropriately.

Here is a pilot  result, which
will be gradually generalized to a very large extent
in the subsequent sections. The first generalization
is given by  Theorem \ref{theorem 8.26.1}
in which the global Lipschitz condition of
$H(z,x)$ with respect to $z$ is replaced with the local one.

\begin{theorem}
                                      \label{theorem 9.12.1}

Under the above assumptions
equation \eqref{9.23.02}
with  boundary
condition $v=g$ on $\partial \Omega$ has a unique
  solution $v\in
C(\bar{\Omega} )\cap C^{1,1}_{ \loc }(\Omega)$.
Furthermore, estimates \eqref{9.22.3}
hold  with $N $ depending only on $\Omega$,
  $K_{0}$, $\Lambda$, and $\delta$
\(in particular, $N $ is independent
of $N'$ and $\rho_{0}$\).
\end{theorem}

\begin{remark}
The constant $N$ in \eqref{9.22.3} is independent of $N'$
 and $\rho_{0}$
which enter  Assumptions \ref{assumption 8.18.2} (i), (iv),
and (v). In   Theorem \ref{theorem 8.29.1}
all these assumptions are dropped. Therefore,
it is worth explaining how we use them here.
Assumption \ref{assumption 8.18.2} (i) will allow us to use
the results about finite-difference equations
while replacing pure second-order derivatives
with  second-order differences and 
Assumptions \ref{assumption 8.18.2}   (iv),
  (v) will allow us to show that the collection of solutions
of finite-difference equations  is almost equicontinuous
in $\Omega$ in the sense specified in 
Lemma \ref{lemma 8.16.1}.

\end{remark}

To prove Theorem \ref{theorem 9.12.1} we need a lemma.
But first of all  we observe that, by standard arguments,
if we have two solutions $u$ and $v$,
then in $\Omega$ (a.e.)
$$
0=\max(H[u],P[u]-K)-\max(H[v],P[v]-K)
$$
$$
=
a^{ij}D_{ij}(u-v)+b^{i}D_{i}(u-v)-c(u-v) 
$$
(cf.~Section \ref{section 3.19.1}), where $a=(a^{ij})$, $b=(b^{i})$, and $c$ are certain 
functions satisfying $a\in\bS_{\delta/2}$, $|b|\leq N'$,
$0\leq c\leq N'$. By the Aleksandrov maximum principle
$u-v=0$ and hence uniqueness.

As is explained above, we are going to use finite 
differences in the proof of the existence and the estimates.
Observe that Assumptions 
\ref{assumption 7.31.1},
\ref{assumption 8.10.1},
\ref{assumption 9.23.01}, and
\ref{assumption 8.8.1}
are satisfied due to 
 Assumptions 
\ref{assumption 8.18.1}, \ref{assumption 8.18.3},  and
\ref{assumption 8.18.2}. 
Assumption \ref{assumption 8.1.1}  is also satisfied 
in light of Assumption \ref{assumption 8.18.2}, in particular, the fact that $H(0,x)$ is a bounded function
follows from Assumption \ref{assumption 8.18.2} (iii).
 
Thus, all the assumptions of Section
\ref{section 9.21.1} are satisfied and we can use
the information obtained in this section
about functions $v_{h}$ defined in $\bar\Omega$
for sufficiently small $h>0$ as unique bounded solutions
of \eqref{2.25.3}-\eqref{2.25.4}.

Estimates \eqref{9.18.9} are estimates of $v_{h}$
on the translates of a multiple of $\bZ^{d}$. They
do not tell us anything about the difference $v_{h}(x)-v_{h}(y)$
if $x$ and $y$ cannot be connected by a broken line
consisting of translates of $hl_{j}$, $|j|=1,...,m$.
However, this information is necessary if we want to
let $h\downarrow 0$ and extract a subsequence
of $v_{h}$ converging everywhere in $\Omega$.
In the proof of the following result
Assumptions \ref{assumption 8.18.2} (iv) and (v)
are crucial.

\begin{lemma}
                                        \label{lemma 8.16.1}
There exists a
constant $N$ such that
for all sufficiently small $h>0$ we have
\begin{equation}
                                           \label{8.16.1}
|v_{h}(x)-v_{h}(y)|\leq N(|x-y|+h)
\end{equation}
whenever $x,y\in \bar\Omega$.

\end{lemma}

Proof. First we want to estimate $v_{h}(x)-v_{h}(y)$
when $|x-y|\leq h$.
Take $x_{0}\in\Omega$ and $h>0$ such that
$4h<\rho_{\Omega}(x_{0}) $. Denote by $N_{1}$ the  
 right-hand side  of \eqref{9.18.9}
and    introduce      
$$
D_{H}=\{(z,x):|z'_{k}|\leq N_{1},|z''_{j}|\leq N_{1} (\rho_{\Omega}(x_{0})-
3h) ^{-1},
$$
$$
|k|=0,...,m,|j|=1,...,m\}\times \bar B_{h}(x_{0}).
$$
We claim that, on the account of our assumptions,
for all sufficiently small $h>0$,
there is a constant $N_{2}$, independent of  $x_{0}$,    
such that,  at all points of differentiability
of $H$ with respect to $(z,x)$ belonging to
the interior of $D_{H}$,
we have 
$$
 \delta/2
\leq D_{z''_{j}}H(z,x)\leq 2\delta^{-1},\quad
|j|=1,...,m,
\quad 0\leq D_{z'_{0}} H(z,x)\leq N_{2},
$$
\begin{equation}
                                         \label{8.17.4}
|D_{(z,x)}H(z,x)|\leq N_{2}.
\end{equation}

Indeed, 
the first two sets of inequalities are given by assumption
and the  Lipschitz continuity of $H(z,x)$ in $z$ uniform with respect to $x$.
If $\rho(x_{0})\leq\rho_{0} /2$ 
($\rho_{0}$ is taken from Assumption \ref{assumption 8.18.2} (v)),
 then $B_{h}(x_{0})
\subset \partial_{\rho_{0}}\Omega$ and $D_{(z,x)}H(z,x)$
(on $D_{H}$)
coincides with the   gradient (in $\bS$) of the right-hand side  
of \eqref{8.17.02}, which is bounded due to
Assumption  \ref{assumption 8.18.2} (i),
 implying
the last inequality in \eqref{8.17.4}.
 However, if $\rho_{\Omega}(x_{0})\geq\rho_{0} /2$,
then (for small $h$) $\rho_{\Omega}(x_{0})-
3h \geq (1/2)\rho_{\Omega}(x_{0})\geq (1/4)\rho_{0}  $ and 
the last inequality in \eqref{8.17.4}
follows from the local Lipschitz continuity of $H(z,x)$
with respect to $(z,x)$,
which in turn is due to the uniform Lipschitz continuity
with respect to $z$ and \eqref{9.31.002}.

The same claim holds true, obviously,  for $P(z'')-K$
and  for
$H(z,x)\wedge(P(z'')-K)$ as well.

Next, it is convenient to continue $v_{h}$ outside $\Omega$
by setting it equal to $g$ and for  a fixed
$l\in\bR^{d}$, such that $|l|\leq h$,
introduce the function $w_{h}(x)=v_{h}(x+l)$. 
Observe that, owing to Lemma   \ref{lemma 9.22.3},  
both points
$$
(v_{h}(x_{0}),\delta_{h,l_{\pm1}}v_{h}(x_{0}),...
,\delta_{h,l_{\pm m}}v_{h}(x_{0}), 
\Delta_{h,l_{\pm1}}v_{h}(x_{0}),...
,\Delta_{h,l_{\pm m}}v_{h}(x_{0}), x_{0})
$$
and
$$
(w_{h}(x_{0}),\delta_{h,l_{\pm1}}w_{h}(x_{0}),...
,\delta_{h,l_{\pm m}}w_{h}(x_{0}), 
\Delta_{h,l_{\pm1}}w_{h}(x_{0}),...
,\Delta_{h,l_{\pm m}}w_{h}(x_{0}), x_{0}+l)
$$
are in $D_{H}$. Furthermore, $x_{0}\in \Omega^{2h}$
and
$$
 H(w_{h}(x_{0}),\delta_{h}w_{h}(x_{0}),
\delta_{h}^{2}w_{h}(x_{0}), x_{0}+l)
\vee(P_{h}[w_{h}](x_{0})-K)=0.
$$
Also  
$$
 H(v_{h}(x_{0}),\delta_{h}v_{h}(x_{0}),
\delta_{h}^{2}v_{h}(x_{0}), x_{0} )
\vee(P_{h}[v_{h}](x_{0})-K)=0.
$$
We subtract these two equations and, by using 
the arbitrariness of $x_{0}$  and what was said
above in the proof, conclude that the function
$u_{h}=w_{h}-v_{h}$ satisfies in $\Omega^{2h}$
$$
a_{k}\Delta_{h,l_{k}}u_{h} +
b_{k}\delta_{h,l_{k}}u_{h}-cu_{h}+f =0 
$$
(cf.~Section \ref{section 3.19.1}),
where $\delta/2\leq a_{k}\leq 2\delta^{-1}$,
 $
|b_{k} |$,  $c\leq N_{2}$, $c\geq0$,  $|f|\leq N_{2}h$.

Next, in Lemma \ref{lemma 8.3.1} reduce $\delta$ if needed,
in such a way that the new one, say $\delta_{1}$,
 will satisfy $N_{2} \leq\delta^{-1}_{1}$ 
and $\delta_{1}\leq\delta/2$,
and then take $\Phi_{0}$ from that lemma.
Then by the comparison principle
we have in $\bR^{d}$ that
$$
|w_{h}-v_{h}|\leq N_{2}h\sup_{\Omega}\Phi_{0}
+\sup_{\bR^{d}\setminus\Omega^{2h}}|w_{h}-v_{h}|.
$$
The last supremum is easily estimated by 
using \eqref{9.13.1} and the fact that
$$
|w_{h}(x)-v_{h}(x)|\leq|v_{h}(x+l)-g(x+l)|
+|g(x+l)-g(x)|+|g(x)-v_{h}(x)|.
$$
This yields \eqref{8.16.1} for $x,y\in\bR^{d}$,
such that $|x-y|\leq h$.

It only remains to observe that,
if $|x-y|\geq h$, one can split the straight
segment between
$x$ and $y$ into adjacent pieces of length $h $ combined with
a remaining one of length less than $h $
and then apply the above result to each
piece (recall that $v_{h}=g$ outside 
$\Omega$). The lemma is proved.  

After that our  theorem is proved in exactly the same way
as
Theorem 8.7 of \cite{Kr11} on the basis  
of Lemma  \ref{lemma 9.22.3}
and the fact that the derivatives of $v$ are weak limits
of finite differences of $v_{h}$ as $h\downarrow0$
(see the proof of Theorem 8.7 of \cite{Kr11}).
One also uses the fact 
that there are sufficiently many pure
 second order derivatives in the directions
of the $l_{i}$'s  to conclude from their boundedness
that the Hessian of $v$ is bounded.     

The following theorem will be used when, after rewriting in terms   
of pure second-order derivatives general
fully nonlinear elliptic
operators, even depending in the Lipschitz continuous way
on the unknown function and its derivatives, we obtain
  an operator that is only locally Lipschitz 
as a function of the unknown function and its first-order derivatives.

\begin{theorem}
                                                \label{theorem 8.26.1}
In Assumption \ref{assumption 8.18.2} \(i\,\)
replace 
the stipulation that $H(z,x)$ ``is Lipschitz continuous with respect to
$z\in \bR^{4m+1}$ with Lipschitz constant $N'$ for any $x\in\bR^{d}$"
with $H(z,x)$  ``is   Lipschitz continuous with respect to $z$
in any ball lying in $\bR^{4m+1}$ with constant independent of $x$''
 and suppose that so
 modified Assumption \ref{assumption 8.18.2}
  along with Assumptions \ref{assumption 8.18.1} 
and \ref{assumption 8.18.3} are satisfied. Then
equation \eqref{9.23.02}
with  boundary
condition $v=g$ on $\partial \Omega$ has a unique
  solution $v\in
C(\bar{\Omega} )\cap C^{1,1}_{ \loc }(\Omega)$.
Furthermore, estimates \eqref{9.22.3}
hold  with $N$  depending only on $\Omega$,
  $K_{0}$, $\Lambda$, and $\delta$.
\end{theorem}

  Proof. Take $N$, that exits by Theorem
\ref{theorem 9.12.1}, and denote 
  by $N^{0}$ the right-hand side 
  \eqref{9.22.3}. Then
take a smooth odd increasing function $\chi(t),t\in\bR$, such that 
$\chi'\leq1$, $\chi(t)
=t$ for $|t|\leq N^{0}$ and $\chi(t)=2N^{0}\sign t$ for 
$|t|\geq 3N^{0}$ and, for $\rho\in(0,\rho_{0}]$,
  introduce
$$
H^{\rho}(z,x)=\delta\sum_{|k|=1}^{m}[z''_{k}-\rho^{-1} \chi
(\rho z''_{k})]
+
H(\chi(z'),\rho^{-1} \chi
(\rho z''  ),x) 
$$
where    
$$
\chi(z')=(\chi(z_{0}'),
\chi(z'_{\pm1}),...
,\chi(z'_{\pm m})),
$$
$$
\chi(\rho z''  )=(\chi(\rho z_{\pm 1}''  ),
  ...
,\chi(\rho z''_{\pm m} )).
$$

Obviously, $H^{\rho}$ is Lipschitz continuous on $\bR^{4m+1}\times\bR^{d}$.
 Also, due to $0\leq\chi'\leq1$, it is easy to check
that Assumption \ref{assumption 8.18.2} (ii\,)
is satisfied for $H^{\rho}$ with the same $\delta$.
 Assumption \ref{assumption 8.18.2} (iii\,) is obviously satisfied 
for $H^{\rho}$ with the same
$K_{0}$ and perhaps smaller $\bar H$.
Assumption  \ref{assumption 8.18.2} (v\,)
is satisfied for $H^{\rho}$ with the same $\rho_{0}$
and
$$
\delta\sum_{|k|=1}^{m}[z''_{k}-\rho^{-1} \chi
(\rho z''_{k} )]
+
H( \rho^{-1} \chi
(\rho z''_{k} ) ) 
$$
in place of $H(z'')$. Finally, 
Assumption  \ref{assumption 8.18.2} (iv)
is satisfied for $H^{\rho}$ with the same $N'$
since $|\rho^{-1} \chi
(\rho z''_{k} )|\leq |z''_{k}|$.

By Theorem \ref{theorem 9.12.1}, there is a solution  
$v^{\rho}$ of the equation
$$
H^{\rho}[v]\vee(P[v]-K)=0
$$
in $\Omega$ with boundary data $v=g$ on $\partial\Omega$
enjoying the properties listed in 
Theorem \ref{theorem 9.12.1}. By the choice of $N^{0}$
and  
\eqref{9.22.3} we have $H^{\rho}[v^{\rho}]=
H [v^{\rho}]$ in $\Omega^{\rho}$ so that
 $v^{\rho}$ also satisfies \eqref{9.23.02}
in $\Omega^{\rho}$. Now we find a sequence $\rho_{n}\downarrow0$
such that $v^{\rho_{n}}$ converge uniformly in $\bar\Omega$
to a function $v$. Obviously, 
$v=g$ on $\partial \Omega$,
  $v\in
C(\bar{\Omega} )\cap C^{1,1}_{ \loc }(\Omega)$
and estimates  \eqref{9.22.3} hold. By fixing $m$
and concentrating on $\Omega^{\rho_{m}}$ by Theorems
3.5.15 and 3.5.9 of \cite{Kr85} we get that 
$v$ satisfies \eqref{9.23.02} in $\Omega^{\rho_{m}}$
for any $m$ and hence in $\Omega$.  This takes care
of the existence and the estimates.

To prove uniqueness it suffices to apply the argument
in Remark \ref{remark 3.10.1} to $\Omega^{\varepsilon}$
for small $\varepsilon>0$
and conclude that in $\Omega$
$$
|u-v|\leq\max_{\partial\Omega^{\varepsilon}}|u-v|,
$$
which, after sending $\varepsilon\downarrow0$
and taking into account that $u,v$ are continuous in $\bar\Omega$
and $u=v=g$ on $\partial\Omega$, yields $u=v$ in $\bar\Omega$.
The theorem is proved.  

\mysection{General elliptic equations with
Lipschitz continuous $H$}

                                         \label{section 8.25.1}
In this section we consider equations not necessarily
written in term of pure derivatives.
Fix some constants $K_{1},K, N'\in[0,\infty)$.
Suppose that we are given a function  $H(u, x)$,
$$
 u=(u',u''),\quad
u'=(u'_{0},u'_{1},...,u'_{d}) \in\bR^{d+1},\quad u''\in\bS,\quad
x\in\bR^{d} . 
$$

\begin{assumption}
                                    \label{assumption 9.23.1}
(i) 
The function   $H(u, x)$ is Lipschitz continuous
with respect to $u''$  and at all points of differentiability  
of $H$ with respect to $u''$ we have
$D_{u''}H \in \bS_{\delta/2}$.

(ii) The number
\begin{equation}
                                                            \label{8.29.9}
\bar{H} :=\sup_{u', x}\big(|H (u',0, x)|-K_{0}|u'|\big)
\quad(\geq0)
\end{equation}
is finite.

(iii) The function   $H(u, x)$ is Lipschitz continuous
with respect to $u'$ and 
at all point of differentiability of $H(u,x)$ with respect to
$u'$ we have
\begin{equation}
                                                     \label{8.28.1}
|D_{u'}H(u,x)|\leq N',\quad D_{u'_{0}}H(u,x)\leq0.
\end{equation}

 (iv)  The function   $H(u, x)$ is locally Lipschitz continuous
with respect to $x$ and 
at all point of differentiability of $H(u,x)$ with respect to
$x$ we have
\begin{equation}
                                    \label{9.31.2}
|D_{x}H (u,x)|\leq N'(1+|u |) .
\end{equation}

(v) For every $u',x$ there exists a subset of $\bS$ of full measure
at every point of which $H(u',u'',x)$ is differentiable
with respect to $u''$ and  
\begin{equation}
                                              \label{9.10.2}
|H  (u, x)-u''_{ij}D_{ u''_{ij}}H(u, x)  |\leq
  K_{1}(\bar H+K+|u'|).
\end{equation} 
 
  (vi) There is a constant $\rho_{0} >0$ and a function $H(u )$ such that
\begin{equation}
                                            \label{8.17.2}
H(u,x)=H(u )
\end{equation}
for all $u$ if $\rho(x)\leq\rho_{0} $.

\end{assumption}

\begin{remark}
                                                   \label{remark 8.30.1}
The condition $D_{u''}H\in \bS_{\delta/2}$ may look strange.
Why do not use $D_{u''}H\in \bS_{\delta }$, which is indeed  imposed
in Theorem \ref{theorem 8.29.1}? The point is that we are going to use the same
$P(u'')$ in Theorems \ref{theorem 9.12.01} and \ref{theorem 8.29.1},
but in the proof of the latter we will replace $H$,
satisfying $D_{u''}H\in \bS_{\delta }$ with the one
satisfying $D_{u''}H\in \bS_{\delta/2}$.
\end{remark}

Here is the   result we are after.  

\begin{theorem}
                                      \label{theorem 9.12.01}

Under the above assumptions
equation \eqref{9.23.02}
  with  boundary
condition $v=g$ on $\partial \Omega$ has a unique
  solution $v\in
C(\bar{\Omega} )\cap C^{1,1}_{ \loc }(\Omega)$.
Furthermore, estimates \eqref{9.22.3} hold
with   
$N$ depending only on $\Omega$,
  $K_{0}$, $K_{1}$, $\Lambda$, and $\delta$
\(in particular, $N$ is independent
of $N'$ and $\rho_{0}$\).

\end{theorem}

\begin{remark}
                                         \label{remark 8.25.1}
If we already know that
equation \eqref{9.23.02}
in $\Omega $  
with  boundary
condition $v=g$ has a solution $v\in
C(\bar{\Omega} )\cap C^{1,1}_{ \loc }(\Omega)$,
then the estimate \eqref{10.15.30}:
\begin{equation}
                                                \label{8.25.1}
\sup_{\bar\Omega}|v|\leq N(\bar{ H}+\|g\|_{C(\Omega )}),
\end{equation}
where $N$ depends only on $\Lambda$, $K_{0}$, $\delta$, and $\diam(\Omega)$,
is obtained on the basis of the maximum principle,
by repeating the proofs of Lemmas \ref{lemma 8.8.2} and
 \ref{lemma 10.6.1} with obvious changes.

\end{remark}

\begin{remark}
                                         \label{remark 10.16.1}

The function $H(  u,x)$ is locally Lipschitz continuous
with respect to $(u,x)$, since, due to  
 Assumption \ref{assumption 9.23.1}, for any $x,y\in\bR^{d}$,
$u,v\in \bR^{d+1}\times \bS$,
$$
|H(u,x)-H(v,y)|\leq |H(u,x)-H(v,x)|+N'(1+|v|)|x-y|
$$
$$
\leq |H(u,x)-H(u',v'',x)|+N'|u'-v'|+N'(1+|v|)|x-y|
$$
$$
\leq N(\delta,d)|u''-v''|+N'|u'-v'|+N'(1+|v|)|x-y|.
$$

\end{remark}

We first prove a version of Theorem \ref{theorem 9.12.01}.

\begin{lemma}
                                           \label{lemma 1.19.1}
In Assumption \ref{assumption 9.23.1} \(v\)
replace \eqref{9.10.2} with
\begin{equation}
                                              \label{2.18.1}
|H  (u, x)-u''_{ij}D_{ u''_{ij}}H(u, x)  |\leq
  K_{1}(\bar H+K+N_{0}+[u']),
\end{equation} 
where
$$
[u']=|(u'_{1},...,u'_{d})|
$$ 
and $N_{0}$ is
the right-hand side of \eqref{8.25.1}. Then the assertion
of Theorem \ref{theorem 9.12.01} holds true.
\end{lemma}

Proof. Without loss of generality we may assume that $K_{1}>0$
(cf.~\eqref{9.10.2}) and define
$$
I=[-K _{1} ,K _{1}],\quad J=[-2K _{1} ,2K _{1}],
\quad C''=I\times\bS_{\delta/2},\quad
B''=J\times
\bS_{\delta/4},
$$
and also recall that $H_{u''}\in\bS_{\delta/2}$.
Then for $u'\in\bR^{d+1}$ and $y''\in\bS$ introduce
$$
B(u',y'', x)=\{(f,l'')\in B'':(\bar{H} +K+N_{0} +[u'])
f+l_{ij}y''_{ij}\leq H(u',y'', x)\}.
$$
Next,  recall  \eqref{4.8.01} and for $u' \in\bR^{d+1}$,
 $x\in\bR^{d}$,
and 
$$
z''=(z''_{\pm1},...,z''_{\pm m})\in \bR^{2m}
$$  ($m$ is the same as in  \eqref{4.8.01}) define
$$
\cH(u',z'',x)
=\inf_{y''\in\bS}\max_{(f,l'')\in
B(u',y'', x)}\big[(\bar{H} +K +N_{0}+[u'])
f
$$
$$
+\sum_{|k|=1}^{m}\lambda_{k}(l'')z_{k}''\big].
$$ 

Then by repeating   word for word the beginning of
Section 4 of \cite{Kr13.1}, we convince ourselves that
the function
$\cH$ is  measurable, Lipschitz continuous
with respect to $ z'' $ with constant
 independent 
of $(u', x)$,
\begin{equation}
                                    \label{9.8.1}
H(u, x)=\cH (u',\langle u''l_{\pm1},l_{\pm1}\rangle,
...,\langle u''l_{\pm m},l_{\pm m}\rangle,  x)
\end{equation}
for all values of the arguments, 
  where $l_{k}$ are taken from
\eqref{6.1.2}, and at all points
of differentiability of $\cH $ with respect to $z''$ we have  
\begin{equation}
                                                \label{9.8.2}
D_{z''}\cH (u',z'', x)
\in [ \hat\delta,\hat\delta ^{-1}]^{2m},
\end{equation}
\begin{equation}
                                      \label{9.8.3}
 (\bar{H} +K+N_{0} +[u'])^{-1}
[\cH (u',z'', x)-\langle z'', 
 D _{ z'' } \cH (u',z'', x)\rangle] \in J
\end{equation}
($\hat\delta$ is introduced in Theorem \ref{theorem 9.26.3}).
Furthermore, $\cH$ is locally Lipschitz continuous with respect
to $(x,u')$ and at all points of its differentiability with
respect to $(x,u')$ we have
\begin{equation}
                                       \label{9.8.4} 
|\cH_{x} (u',z'', x) |
\leq N (1+|z''|+|u'|),
\end{equation}
\begin{equation}
                                       \label{2.19.1}
|\cH_{u'}(u',z'', x)|\leq N(1+|u'|+|z''|),
\end{equation}
where $N$ is a constant independent of $u',z'', x  $
(for getting \eqref{2.19.1} assuming the first
inequality in \eqref{8.28.1}
is crucial).

Finally (what is not coming from Section 4 of \cite{Kr13.1}),
the function $\cH$ is a decreasing function of $u'_{0}$.
This follows from the fact that, since $[u']$
is independent of $u'_{0}$ and $H(u',y'',x)$
is a decreasing function of $u'_{0}$, for smaller
values of $u'_{0}$ the set $B(u',y'', x)$
is smaller.

Now we want to apply Theorem \ref{theorem 8.26.1}.
Since by assumption $e_{i}\in\Lambda$, $i=1,...,d$,
there is an identification of some $e_{i}$ with some $l_{k}$.
Let $e_{i}=l_{k(i)}$, for $i=1,...,d$,
and then   for 
$$
z=(z',z''),\quad z'=(z'_{0}, z'_{\pm1},
...,z'_{\pm m})\in\bR^{2m+1},\quad z''=(z''_{\pm1},...,
z''_{\pm m})\in\bR^{2m}
$$
and $x\in\bR^{d}$, set
$$
\hat \cH(z,x)=
\cH (z_{0}',z'_{k(1)},...,z'_{k(d)},z'',x).
$$

Observe that, if $\rho_{\Omega}(x)\leq\rho_{0}$,
$\hat H(z,x)$ depends only on $z $, since
then $H(u,x)$ depends only on $u $ by Assumption
\ref{assumption 9.23.1} (vi). Assumption \ref{assumption 8.18.2}
(iv) is satisfied for $\hat H(z,x)$ due to
  \eqref{2.19.1}. 
If we take $z''=0$ in \eqref{9.8.3},
then we see that
$$
|\hat H(z',0,x)|\leq 2
K_{1}(\bar H+K +N_{0}+|z'|).
$$
Hence, Assumption \ref{assumption 8.18.2}
(iii) is satisfied for $\hat H(z,x)$
with
$$
2K_{1}(\bar H+K +N_{0})\quad\text{and}\quad 2K_{1}
$$
in place of $\bar H$ and $K_{0}$. 
Assumption \ref{assumption 8.18.2}
(ii\,) is satisfied for $\hat H(z,x)$
with $\hat\delta$ in place of $\delta$ owing to \eqref{9.8.2}.
Finally, the modification of
Assumption \ref{assumption 8.18.2} (i\,) stated in
Theorem \ref{theorem 8.26.1} is satisfied
for $\hat H(z,x)$ in light of \eqref{9.8.2},
 \eqref{9.8.4}, and \eqref{2.19.1}. Hence, all the assumptions
of Theorem \ref{theorem 8.26.1} are satisfied.
We draw the reader's attention to the fact that
the assumptions
of Theorem \ref{theorem 8.26.1} are satisfied with $\hat\delta$
in place of $\delta$ (among a few other substitutions) and that
$P[u]$ 
in Theorem  \ref{theorem 8.26.1} is taken from
\eqref{6.1.1}, in accordance with which $P[u]$ in Theorems
 \ref{theorem 8.29.1} and \ref{theorem 9.12.01}
is of type \eqref{6.1.1} with $\hat\delta$ in place of~$\delta$.

By Theorem \ref{theorem 8.26.1} the equation
\begin{equation}
                                                     \label{8.29.5}
\hat H(v(x),D_{l_{k(i)}}v(x),D^{2}_{l_{k}}v(x),x)
\vee (P[v](x)-K)=0
\end{equation}
in $\Omega$ with boundary condition $v=g$
on $\partial\Omega$ has a unique solution
$v\in C(\bar\Omega)\cap C^{1,1}_{ \loc }(\Omega)$ and the   
estimates \eqref{9.22.3}
 hold true if the right-hand side is
replaced with
 \begin{equation}
                                                         \label{8.29.6}
N[K_{1}(\bar H+K +N_{0})+K+\|g\|_{C^{1,1}(\bR^{d})}],
\end{equation}
where $N$ is a constant depending only on $\Omega$,
  $K_{1}$, $\Lambda$, and $\hat \delta$. We recall what $N_{0}$ is
and estimate \eqref{8.29.6} from above by 
$$
N( \bar{H}+K 
+\|g\|_{ C^{1,1} (\bR^{d})})
$$
with $N$ as in the statement of the theorem.
Then it only remains to observe that, in light of \eqref{9.8.1},
equation \eqref{8.29.5} coincides with \eqref{9.23.02}.
This proves the lemma.

{\bf Proof of Theorem \ref{theorem 9.12.01}}. 
Introduce $N_{0}$ as the right-hand side of \eqref{8.25.1}  
and cut-off the range of the variable $u'_{0}$ by setting
$$
H^{0}(u,x)=H(-N_{0}\vee u'_{0}\wedge N_{0},u'_{1},...,u'_{d},
u'',x).
$$
Then, owing to \eqref{9.10.2},  
$$
|H^{0}(u,x)-u''_{ij}D_{u''_{ij}}H^{0}(u,x)|
\leq K_{1}(\bar H+K +N_{0}+[u']).
$$
 
By Lemma \ref{lemma 1.19.1} equation
\begin{equation}
                                              \label{2.20.1}
H^{0}[v]\vee(P[v]-K)=0
\end{equation}
in $\Omega$ with boundary condition $v=g$
on $\partial \Omega$ has a unique solution
as stated in Theorem \ref{theorem 9.12.01}.
Furthermore, since
$$
|H^{0}(u',0,x)|\leq \bar H+K_{0}|u'|
$$
and $H^{0}$ is a decreasing function of $u'_{0}$,
by Remark \ref{remark 8.25.1} we have $|v|\leq N_{0}$.
But then \eqref{2.20.1} coincides with \eqref{9.23.02}
and the theorem is proved.

\mysection{Proof of Theorem \protect\ref{theorem 8.29.1}}
                                          \label{section 8.29.1}

\begin{remark}
                                                 \label{remark 11.3.1}
If we already know that there is a
$v\in
C(\bar{\Omega} )\cap C^{1,1}_{ \loc }(\Omega)$ satisfying
equation \eqref{9.23.02} and such  that $v=g$ on $\partial\Omega$,
 then
estimate   \eqref{10.18.3} with $N$ depending only
on $d,\delta, K_{0}$ and $\diam(\Omega)$ is  
obtained by an easy adaptation of the proofs of Lemmas \ref{lemma 8.8.2} and
 \ref{lemma 10.6.1}.

\end{remark}

\begin{remark}
                                           \label{remark 2.6.1}
Assertion (ii) follows from (i) and  
  Theorem 1.2 of \cite{DKL}.  To show this, 
observe that (cf. the proof of Lemma \ref{lemma 9.3.10})
$$
H(u,x)=[H(u,x)-H(u',0,x)]+H(u',0,x)=a^{ij}u''_{ij}
+\sum_{i=0}^{d}b^{i}u'_{i}+\theta\bar H,
$$
where $\bS_{\delta}$-valued $a=(a^{ij}) $, $\bR^{d+1}$-valued
 $b=(b^{i})$, and $[-1,1]$-valued $\theta$ are certain functions
of $(u,x)$ such that $|b|\leq K_{0}$. It follows from the construction
of $P$ that for a constant $\kappa=\kappa(\delta,d)>0$ we have
$$
H(u,x)\leq P(u)-\kappa \sum_{i,j}|u''_{ij}| +\bar H+K_{0}|u'|.
$$
Hence, if we have a solution $v$ like in assertion (i),
then
\begin{equation}
                                                \label{3.6.5}  
H[v]\leq P[v]-\kappa \sum_{i,j}|D_{ij}v|+N_{0},
\end{equation}
where $N_{0}$ is $\bar H$ plus $K_{0}$
times the right-hand side of \eqref{8.29.10}.

Next,
$$
\max(H[v],P[v]-K)=P[v]+Q[v],
$$
where $Q[v]= (H[v]-P[v]+K)_{+}-K$ and, owing to  
 \eqref{3.6.5},
$Q[v]=-K$ if
$$
-\kappa \sum_{i,j}|D_{ij}v|+N_{0}\leq-K,\quad
\kappa \sum_{i,j}|D_{ij}v|\geq N_{0}+K.
$$
If the opposite inequality holds, then
\begin{equation}
                                            \label{1.14.4}
|Q[v]|\leq | H[u]-H(v,Dv ,0,x)|+|P[v]|+N_{0}+2K
\leq NN_{0},
\end{equation}
where $N$ depends only on $\delta$ and $d$.
It follows that the inequality between
the extreme terms in \eqref{1.14.4} holds 
and $v$ is a $W^{2}_{p,\loc}(\Omega)\cap C(\bar \Omega)$-solution
of $P[u]=-f$, where $f=-Q[v]$ is in $L_{p}(\Omega)$
(actually, bounded).
By Theorem 1.2 of \cite{DKL} this equation with boundary data
$g$ has a  solution   $u\in W^{2}_{p }(\Omega) $
and
$$
\|u\|_{W^{2}_{p}(\Omega)}
\leq N\|f\|_{L_{p}(\Omega)}.
$$
By uniqueness, $v=u$, and we get \eqref{1.13.2}.

\end{remark}

Because of Remark \ref{remark 2.6.1} below we are only
dealing with assertion (i) of  
Theorem \ref{theorem 8.29.1}.
The proof of it is based on
Theorem \ref{theorem 9.12.01} and will be achieved
in several steps in the first two of which  we drop
Assumptions \ref{assumption 9.23.1}   (v), (vi\,),
one by one.

\begin{lemma}
                                                      \label{lemma 8.29.4}
In addition to the the assumptions of Theorem \ref{theorem 8.29.1}
let the assumptions of Theorem \ref{theorem 9.12.01},
apart from Assumption \ref{assumption 9.23.1} \(vi\,\),
be satisfied. Then the assertions of Theorem 
\ref{theorem 9.12.01} are still true.
\end{lemma}

Proof. Define $H(u'')=H(0,u'',0)$ and for $\rho_{0}>0$
find a function $\zeta_{\rho_{0}}\in C^{\infty}_{0}(\Omega)$
such that $\zeta_{\rho_{0}}=0$ on $\Omega\setminus\Omega^{\rho_{0}}$  
and $\zeta_{\rho_{0}}=1$ on $\Omega^{2\rho_{0}}$ and $0\leq
\zeta_{\rho_{0}}\leq1$ in $\Omega$.
Then introduce
$$
H^{\rho_{0}}(u,x)=\zeta_{\rho_{0}} H(u,x)+(1-\zeta_{\rho_{0}})H(u'').
$$

Obviously, $H^{\rho_{0}}(u,x)=H(u'')$ if $\rho(x)\leq\rho_{0}$
and condition \eqref{9.10.2} is certainly satisfied for
$H^{\rho_{0}}$ with $2K_{1}$ in place of $K_{1}$.
Furthermore, 
$$
|H(u,x)|\leq|H(u,x)-H(0,u'',x)|+|H(0,u'',x)-H(0,x)|
$$
$$
+|H(0,x)|
\leq N'|u'|+N|u''|+\bar H,
$$
where $N'$ is taken from \eqref{8.28.1} and $N$ accounts
for the Lipschitz continuity of $H$ in $u''$.
It follows that condition \eqref{9.31.2}
is satisfied for $H^{\rho_{0}}$ with a different constant $N'$
(by the way, depending on $\rho_{0}$ but this is irrelevant).
Condition \eqref{8.28.1} is satisfied for $H^{\rho_{0}}$ with 
the same $N'$ and the number $\bar H^{\rho_{0}}$ is at most twice
the $\bar H$ from \eqref{8.29.9}.

By Theorem \ref{theorem 9.12.01} the equation
\begin{equation}
                                                        \label{8.29.11}
\max(H^{\rho_{0}}[v^{\rho_{0}}],P[v^{\rho_{0}}]-K)=0
\end{equation}
in $\Omega$ with boundary condition $v^{\rho_{0}}=g$ on $\partial
\Omega$ has a unique solution $v^{\rho_{0}}$
with the properties described in that theorem
(with $N$ from Theorem \ref{theorem 9.12.01} in \eqref{9.22.3}
 multiplied by 2,
because $\bar H^{\rho_{0}}\leq2\bar H$).
One sends $\rho_{0}\downarrow0$ and finishes the proof of existence
of a solution with desired properties as in the end of the proof
of Theorem \ref{theorem 8.26.1}. Uniqueness is also
shown as there.  
The lemma is proved. 

\begin{lemma}
                                             \label{lemma 8.30.1}
In addition to the the assumptions of Theorem \ref{theorem 8.29.1}
let the assumptions of Theorem \ref{theorem 9.12.01},
apart from Assumptions \ref{assumption 9.23.1}   \(v\,\), \(vi\,\),
be satisfied. Then the assertions of Theorem 
\ref{theorem 9.12.01}  are still true
with $N$ in \eqref{9.22.3}
this time independent of $K_{1}$
\(which does not even enter the assumptions of the present lemma\).
\end{lemma}

Proof. Introduce 
$$
 P_{0}(u'')=\max_{a\in\bS_{\delta/2}}
 \tr a  u'' ,\quad H_{K}=\max(H,P_{0}-K)  .
$$ 
Obviously, $P_{0}(u'')
\leq P(u'')$, so that the equation
$$
\max(H_{K}[u],P[u]-K)=0
$$
is equivalent to \eqref{9.23.02}. Furthermore,
as is shown in Section 3 of \cite{Kr13.1} (or as follows
from Section \ref{section 3.19.1}),
the function $H_{K}$ satisfies Assumption
 \ref{assumption 9.23.1} (i) (here Remark \ref{remark 8.30.1}
is relevant).

 It is also shown that the function $H_{K}$   satisfies 
Assumptions
 \ref{assumption 9.23.1} (iii), (iv) with the same constant $N'$
and the number $\bar H_{K}\leq\bar H<\infty$, so that
Assumption
 \ref{assumption 9.23.1} (ii) is also satisfied  for $H_{K}$.

By Lemma 3.2 of  \cite{Kr13.1},
   Assumption
 \ref{assumption 9.23.1} (v) is   satisfied for $H_{K}$
with
$$
K_{1}=N(K_{0}+1) ,
$$
where $N$ depends only on  $d$ and $\delta$.

Now to finish the proof of the lemma it only
remains to refer to Lemma \ref{lemma 8.29.4}.
The lemma is proved.  

\begin{remark}
One might think that we could take 
$$
H_{K}(u,x)=\max(H(u,x),P(u'')-K)
$$
replacing $P_{0}$ with $P$. However, then $H_{K}$
would satisfy  Assumption
 \ref{assumption 9.23.1} (i) with $\hat\delta/2$ in place of $\delta/2$
and this would make  Lemma \ref{lemma 8.29.4} inapplicable.
\end{remark}

We finish the proof of Theorem \ref{theorem 8.29.1}
by  a result, that is   equivalent to this theorem.

\begin{lemma}
                                            \label{lemma 8.31.2}
In addition to the the assumptions of Theorem \ref{theorem 8.29.1}
let the assumptions of Theorem \ref{theorem 9.12.01},
apart from Assumptions \ref{assumption 9.23.1}
\(iii\,\),  \(iv\,\), \(v\,\), \(vi\,\),
be satisfied. Then the assertions of Theorem 
\ref{theorem 9.12.01}  are still true
with $N$ in \eqref{9.22.3}
  independent of $K_{1}$ \(which is nowhere to be found
in the assumptions of the present lemma\).
\end{lemma}

The proof of this lemma is achieved by repeating the proof
of Lemma 3.1 of \cite{Kr13.1} by mollifying $H$ with respect
to $x$ and $u'$ and using elliptic versions
of the theorems about passage to the limit
inside nonlinear operators (see, for instance, Section
3.5 in \cite{Kr85})
 instead
of parabolic ones, that were used in \cite{Kr13.1}.
These mollifications allow one to rely on Lemma
\ref{lemma 8.30.1}. We say more along these lines
in Section \ref{section 3.14.1}.

   \mysection{Proof of Theorems \protect\ref{theorem 2.20.1}
and \protect\ref{theorem 2.20.2}}

                                              \label{section 3.14.1}

{\bf Proof of Theorem \ref{theorem 2.20.1}}.
First assume that $g=0$. We are going to use the
following result which is
  the correct version of  Theorem 5.4 of \cite{Kr13}.
  We would not need the correction if in \cite{Kr13}
and here
we assumed that $F_{u''}\in\bS_{\delta}$. Without this assumption
the claim made in the end of the proof of Theorem 5.3 of \cite{Kr13}, 
on which Theorem 5.4 of \cite{Kr13} is based, that 
$F[u]=F[u]-F[0]=Lu$ for an $L\in\bL_{\delta,K_{0}}$
is unsubstantiated. However,  for $u$ satisfying \eqref{7.29.1}
we have
$$
-H(u,Du,0,x)=H[u]-H(u,Du,0,x)=a_{ij}D_{ij}u,
$$
where $(a_{ij})$ is an $\bS_{\delta}$-valued function and
$|H(u,Du,0,x)|\leq \bar G+K_{0}(|u|+|Du|)$. This is
enough to get the correct version of Theorem 5.3 of \cite{Kr13}
which 
along with 
interpolation theorems 
lead  to the correct version of Theorem 5.4 of \cite{Kr13},
which we present below,
right away.  
Recall that $p>d$ and $\WO^{2}_{p}(\Omega)$ is the subset
of $W^{2}_{p}(\Omega)$ of functions vanishing on $\partial\Omega$.

\begin{theorem}
                                   \label{theorem 11.28.1}
Let $p\in(d,\infty)$. Then there exists
a constant $\theta>0$ depending only on $\Omega$, $p,d$, $K_{F}$,
 and $\delta$
such that if Assumption \ref{assumption 2.13.4}
is satisfied with this $\theta$,
then for any function
 $
u\in \WO^{2}_{p}(\Omega)
 $
 satisfying \eqref{7.29.1} 
we have
\begin{equation}
                                         \label{11.28.2}
\|u\|_{W^{2}_p(\Omega)}\le N
\| F[u]\|_{L_p(\Omega)}+ N\|\bar G\|_{L_p(\Omega)}
+N\|u\|_{L_p(\Omega)}  +Nt_{0},
\end{equation}
where $N$ depends only on  $\Omega,R_{0},d,p $, $K_{F}$, 
and $\delta$. 
 \end{theorem}

Since $|F[u]|=|G[u]|\leq  
\hat\theta|D^{2}u|+ \bar G+K_{0}[|u|+|Du|]$,
using interpolation theorems allowing to estimate
the $L_{p}$-norm of $Du$ through the $L_{p}$-norms of $D^{2}u$
and $u$, we immediately get \eqref{eq16.01}  
after choosing $\hat\theta$ sufficiently small.

 In the general case introduce $\hat g(x)=(g(x),Dg(x),D^{2}g(x))$
and
$$
\hat H(v,x)=H(v+\hat g(x),x) ,\quad w(x)=u(x)-g(x).
$$
Observe that $\hat H[w]=0$ in $\Omega$ (a.e.) and
$w\in \WO^{2}_{p}(\Omega)$. Furthermore, for
$$
\hat G(v,x):=\hat H(v,x)-F(v'',x)
$$
we have
$$
|\hat G(v,x)|=|F(v''+D^{2} g(x),x)-F(v'',x)
+G(v+\hat g(x),x)|
$$
$$
\leq \hat\theta|v''|+ N|D^{2}g(x)|+\bar G+K_{0}
\big(|v'_{0}+g(x)|^{2}
$$
$$
+\sum_{i=1}^{d}|v'_{i}
+D_{i}g(x)|^{2}\big)^{1/2}
\leq\,\hat\theta|v''|+ \hat{\!\bar G}+K_{0}|v'|,
$$
where 
$$
\hat{\!\bar G}=N|D^{2}g |+\bar G+K_{0}
 ( g |^{2}+|Dg|^{2}))^{1/2}
$$
and $N$ depends only on $K_{F}$ and $d$.

It follows that the above result is applicable to $w$
which  leads to \eqref{eq16.01} in the general case. 
 The theorem is proved.

{\bf Proof of Theorem \ref{theorem 2.20.2}}.
First assume that $g\in C^{1,1}(\bR^{d})$
and $\omega(t)=N_{0}t$ in Assumption \ref{assumption 3.11.2} (ii),
where $N_{0}$ is a constant.
  In that case we closely follow a few steps in the proof
of Theorem 2.1 of \cite{Kr13} given in Section 6 there.
 Introduce a function of one variable by setting
$\xi_{K}(t)=0$ for $|t|\leq K$ and
$\xi_{K}(t)=t$ otherwise, set $H^{0}_{K}(x)=\xi_{K}(H(0,x))$
and define
$$
H_{K}(u,x)= H(u,x)-H^{0}_{K}(x).
$$

Notice that since $F(0,x)=0$, we have $|H^{0}_{K}|=
|\xi_{K}(G(0,x))|\leq \xi_{K}(\bar{G})
\leq\bar{G}$.
Also 
$$
H_{K}(0,x) =G(0,x)I_{|G(0,x)|\leq K}
$$
so that $|H_{K}(0,x)|\leq K\wedge\bar G(x)$ and
$$
|H_{K}(u',0,x)|\leq|H_{K}(u',0,x)-H_{K}(0,x)|+
|H_{K}(0,x)|\leq N_{0}|u'|+K.
$$
This shows that, for $H_{K}$, Assumption
\ref{assumption 8.30.1} (ii) is satisfied
with $N_{0}$ in place of $K_{0}$
and an $\bar H\leq K$. Assumption
\ref{assumption 8.30.1} (i) is satisfied
with  the same $\delta$ and $\omega$.
By Theorem \ref{theorem 8.29.1}, with $P(u'')$ from   \eqref{3.28.1}, there is a  
solution $v_{K}\in W^{2}_{p}(\Omega)$ of the equation
\begin{equation}
                                                 \label{3.11.1}
\max(H_{K}[v_{K}], P[v_{K}]-K)=0
\end{equation}
in $\Omega$ with boundary data $g$.

We want to apply Theorem \ref{theorem 2.20.1}
to \eqref{3.11.1}. To this end introduce
$$
\hat H_{K}(u,x)=\max(H_{K}(u,x),P(u'')-K),
$$
$$
F_{K}(u'',x)=\max(F(u'',x),P(u'')-K),\quad G_{K}(u,x)=
\hat H_{K}(u,x)-F_{K}(u'',x).
$$
Below in this section by $N$ we denote
various constants which depend  only on $\Omega,R_{0},d,p $,
$K_{0}$, and $\delta$.

Observe that
$$
|G_{K}(u,x)|\leq |H(u,x)-H^{0}_{K}(x)-F(u'',x)| 
$$
\begin{equation}
                                                    \label{11.30.1}
=
|G(u,x)-H^{0}_{K}(x)|\leq  \hat\theta|u''|+ 
 K_{0}|u'|+2\bar{G}(x).
\end{equation}
Furthermore, $F_{K}$ obviously satisfies Assumption  
\ref{assumption 2.13.4} (i) perhaps with  a constant $N$
(independent of $K$) in place of $K_{F}$.
 To check that the remaining conditions
in Assumption \ref{assumption 2.13.4}  are satisfied take
$z\in\Omega$, $r\in(0,R_{0}]$, the function $\bar{F}=\bar{F}_{z,r}$
from Assumption \ref{assumption 2.13.4}  and set
$$
\bar{F}_{K}(u'')=\max(\bar{F}(u''),P(u'')-K).
$$
Notice that
$$
t^{-1}
|F_{K}(tu'' ,x)-\bar{F}_{K}(tu'')| \leq
t^{-1}
|F(tu'' ,x)-\bar{F}(tu'')|,
$$
which implies that Assumption \ref{assumption 2.13.4}  is satisfied
indeed with the same $\theta$, $R_{0}$, and $t_{0}$ and, perhaps,
modified $\delta>0$ independent of $K$.

It follows by Theorem \ref{theorem 2.20.1}
\begin{equation}
                                                \label{11.28.4}
\|v_{K}\|_{W^{2}_{p}(\Omega)}
\leq N \|\bar{G}\|_{L_{p}(\Omega)}
+N\|g\|_{W^{2}_{p}(\Omega)}
+N\|v_{K}\|_{L_{p}(\Omega)} +N t_{0} .
\end{equation}  
Since equation \eqref{3.11.1} can be rewritten as
$$
a^{ij}D_{ij}v_{K}+b^{i}D_{i}v_{K}-cv_{k}+2\theta\bar G=0
$$
(cf.~Lemma \ref{lemma 9.3.10}), where $(a^{ij})$
is a certain $\bS_{\bar\delta}$-valued function,
$|(b^{i})|\leq K_{0}$, $c\geq0$, by the Aleksandrov estimate
we can eliminate $\|v_{K}\|_{L_{p}(\Omega)}$ on the right 
in \eqref{11.28.4} and conclude that
\begin{equation}
                                                \label{3.13.4}
\|v_{K}\|_{W^{2}_{p}(\Omega)}
\leq N \|\bar{G}\|_{L_{p}(\Omega)}
+N\|g\|_{W^{2}_{p}(\Omega)} +N t_{0} .
\end{equation}

In this way we completed a crucial step consisting of
 obtaining a uniform
control of the $W^{2}_{p}(\Omega)$-norms of $v_{K}$.

We now let $K\to\infty$. As is well known, 
there is a sequence $K_{n}\to\infty$
as $n\to\infty$ and $v\in W^{2}_{p}(\Omega)$
  such that $v_{K}\to v$ weakly in $W^{2}_{p}(\Omega)$.
Of course, estimate \eqref{3.13.4} holds with $v$
in place of $v_{K}$.

By the compactness of embedding of
$W^{2}_{p}(\Omega)$ into $C(\bar{\Omega})$ we have that
  $v_{K}\to v$ also uniformly.

Next, the operator $H[u]$ fits in  the scheme of Section
5.6 of \cite{Kr85} as long as $\omega(t)=N_{0}t$.
 In addition, by recalling that 
$|H^{0}_{K}|\leq \xi_{K}(\bar{G})$ we get
$$
|H[v_{K}]|=|\max(H[v_{K}]-H^{0}_{K},P[v_{K}]-K)-H[v_{K}]|
$$
$$
=|\max(0,P[v_{K}]-H[v_{K}]+H^{0}_{K}-K)-H^{0}_{K}|
$$
$$
\leq(P[v_{K}]-H[v_{K}]+H^{0}_{K}-K)_{+}+|H^{0}_{K}|
$$
$$
\leq(P[v_{K}]-H[v_{K}]+H^{0}_{K} -K)_{+}+ \xi_{K}(\bar{G})
$$
$$
\leq(N|D^{2}v_{K}|+
N|Dv_{K}|+N|v_{K}| +\bar{G}-K)_{+}+ \xi_{K}(\bar{G}),
$$
so that
$$
\| H[v_{K}] \|_{L_{d}(\Omega)}^{d}\leq NK^{d-p}
\int_{\Omega}
( |D^{2}v_{K}|+
 |Dv_{K}|+ |v_{K}|+\bar{G})^{p}\,dx\to0
$$
as $K\to\infty$. By combining all these facts and applying
  Theorems 3.5.15 and 3.5.6
of \cite{Kr85} we conclude that $H[v]=0$ and this finishes
the proof of   the theorem if $\omega=N_{0}t$ and
$g\in C^{1,1}(\bR^{d})$.  If $g=0$,
in this way we obtain a slightly different proof of the existence
part in Theorem 2.1 of \cite{Kr13}.

 To pass to the general $\omega$, take a nonnegative
$\zeta\in C^{\infty}_{0}(\bR^{d+1})$, which integrates
to one and has support in the unit ball centered at the origin
and   define $H^{n}(u,x)$  as the convolution of 
$H (u,x)$ and $n^{ d+1}\zeta(n u')$  performed
with respect to $u'$.
Keep $F_{n}=F$. As is easy to see, for each $n$, $H_{n}$
satisfies Assumption \ref{assumption 10.5.1} with the same $K_{0}$
 and $\hat\theta$ 
and $\bar G+1/n$ in place of $\bar G$. 
Furthermore, for any $k=0,...,d$
$$
H^{n}_{ u'_{k}}(u',u'', x)=n\int_{\bR^{ d+1}}
 H(u'-v'/n,u'', x ) \zeta_{u'_{k}}(v' )\,dv'
$$
$$
=n\int_{\bR^{ d+1}}
[H(u'-v'/n,u'', x )-H(u',u'', x)]
 \zeta_{u'_{k}}(v' )\,dv'.
$$
It follows that
$$
|H^{n}_{ u'_{k}}(u ,t,x)|\leq n\omega(1/n)|B_{1}|_{\bR^{d+1}}\sup| D\zeta|,
$$
where $|B_{1}|_{\bR^{d+1}}$ is the volume of the unit ball in $\bR^{d+1}$,
so that $H^{n}$ also satisfies Assumption \ref{assumption 3.11.2} (ii)
with a constant $N$, depending on $n$, times $t$
in place of $\omega$. Assumptions \ref{assumption 3.11.2} (i)
and \ref{assumption 10.30.1} are obviously satisfied for $H^{n}$.
Since we did not change $F$, Assumption \ref{assumption 2.13.4}
is satisfied and by the first part of the proof,
there exists $v^{n}\in W^{2}_{p}(\Omega)$ such that
$H^{n}[v^{n}]=0$ in $\Omega$ and $v^{n}=g$ on
$\partial\Omega$. Furthermore, \eqref{3.13.4} holds
with   $v^{n}$ in place of $v_{K}$,
$\|\bar G\|_{L_{p}(\Omega)}+1$ in place of
$\|\bar G\|_{L_{p}(\Omega)}$, and $N$ independent of $n$.

Owing to embedding theorems, $W^{2}_{p}(\Omega)
\subset C^{1+\gamma}(\Omega)$, where $\gamma>0$,
the sequence $\{v^{n},Dv^{n}\}$,       
being uniformly bounded and uniformly continuous, has a subsequence uniformly in $\Omega$ converging
to  $v,Dv$, where  $v  \in  
W^{2}_{p}(\Omega)$. For simplicity 
of notation we suppose that the whole sequence $v^{n}$, $Dv^{n}$
converges. Of course, $v=g$ on $\partial\Omega$.

Observe that for $m\geq n$
\begin{equation}
                                              \label{10.3.4}
 \check{H}^{n} [v^{m}]\geq0
\end{equation}
in $\Omega$ (a.e.), where
$$
\check{H}^{n} (u, x):=\sup_{k\geq n}
 H^{k}(v^{k}( x),Dv^{k}( x),u'', x) .
$$

In light of \eqref{10.3.4} and the fact that
 the norms $\|v^{n}\|_{W^{ 2}_{p}(\Omega)}$ are bounded,
by Theorems 3.5.15 and  3.5.9  of \cite{Kr85}
 we have 
\begin{equation}
                                              \label{10.3.2}
 \check{H}^{n} [v ]\geq0 
\end{equation}
in $\Omega $ (a.e.).

Now we notice that  
$$
|H^{k}(u, x)-H (u, x)|\leq\omega(1/k),
$$
$$
|H^{k} (v^{k}( x),Dv^{k}( x),D^{2}v( x), x)
-H(v ( x),Dv ( x),D^{2}v( x), x)|
$$
$$
\leq\omega\big(|v^{k}-v|( x)+|Dv^{k}-Dv|( x)
 \big)+\omega(1/k),
$$
which along with what was said above implies that
\begin{equation}
                                              \label{9.22.5}
 \hat{H}^{n} [v ]\geq-\varepsilon_{n}
\end{equation}
in $\Omega $ (a.e.), where the functions $\varepsilon_{n}
\to 0$ in $\Omega $ (even   uniformly) and
$$
\hat{H}^{n} (u, x):=\sup_{k\geq n}
H^{k}(u, x), 
$$
By letting $n\to\infty$ in \eqref{9.22.5} and using that
$H^{k}(u,x)\to H(u,x)$ as $k\to\infty$ for any $(u,x)$,
we conclude that $H[v]\geq0$ in $\Omega $ (a.e.).

The inequality $H[v]\leq0$ in $\Omega $ (a.e.)
is obtained by similar arguments starting with
$$
 \inf_{k\geq n}
 H^{k}(v^{k}( x),Dv^{k}( x),u'', x) .
$$
 
The passage from $g\in C^{1,1}(\bR^{d})$
to $g\in W^{2}_{p}(\Omega)$ 
is achieved by mollifying $g$
and using a very simplified version of the above
arguments. The theorem is proved.

\mysection{Appendix}

                                            \label{section 3.19.1}

Here are two results used mostly inexplicitly in various combinations
at various places
in the article.

\begin{lemma}
                                                  \label{lemma 7.19.2}
 
Let $F(u )$ be a real-valued    Lipschitz continuous
function defined in $\bR^{n}$. Assume  
 that  there is a
convex closed bounded set $A \subset\bR^{n}$ and a set of full measure
$D'_{F}\subset \bR^{n}$
such that  at all points $u\in D'_{F}$ the function   $F$ is differentiable 
 and $DF (u )\in A $. Then
for any   $u,v\in \bR^{n}$  
there exists an $a\in A$ such that
\begin{equation}
                                                       \label{7.19.2}
F(u )-F(v )=a^{i} (u^{i} -v^{i} ).
\end{equation}
In particular,
\begin{equation}
                                                       \label{7.19.1}
\min_{a\in A}[a^{i} (u^{i} -v^{i} )]\leq
F(u )-F(v )\leq \max_{a\in A}[a^{i} (u^{i} -v^{i} )].
\end{equation}

\end{lemma}

Proof.   Fix $v\in \bR^{n}$.
By using the Fubini theorem in polar coordinates we obtain
that, for almost all points $\omega\in\partial B_{1}$, $v+t\omega\in
D'_{F}$
for almost all $t$. In addition, $F(v+t\omega)$ 
is a Lipschitz 
and absolutely continuous function of $t$ on that interval. It follows that,
for almost all points $\omega\in\partial B_{1}$, for almost all $t$
  we have
$$
 \partial_t  F(v+t\omega)=\omega^{i}[D_{u^{i}}F](v+t\omega)
$$
and for all $t>0$  
$$
F(v+t\omega)=F(v)+ta^{i}\omega^{i},
$$
where
$$
a^{i}= \frac{1}{t}\int_{0}^{t}[D_{u^{i}}F](v+s\omega)\,ds.
$$
Since $A$ is closed and convex, $a\in A$, and 
we obtain \eqref{7.19.2} for almost all $u\in \bR^{n}$. Then by compactness of $A$
and the continuity of $F$ we extend \eqref{7.19.2} to all 
$u,v\in \bR^{n}$.
The lemma is proved.

The following is in a sense a converse 
statement to Lemma \ref{lemma 7.19.2}.

\begin{lemma}
                                          \label{lemma 7.24.1}
Let $A$   be as in Lemma \ref{lemma 7.19.2} and let
$F(u )$ be a real-valued function such that for all $u,v\in \bR^{n}$
\begin{equation}
                                                     \label{7.24.1}
F(u )-F(v )\leq \max_{a\in A}[a^{i} (u^{i} -v^{i} )]=:U(u-v).
\end{equation}
Then
$F$ is Lipschitz continuous on $\bR^{n}$ and at all points $u
\in \bR^{n}$, at which $F$ is differentiable,
 we have $DF (u )\in A $.
\end{lemma}

Proof.   Obviously $|U(w)|\leq N|w|$, where $N$ is a constant
independent of $w$. By interchanging $u$ and $v$
in \eqref{7.24.1} we see that $F(u )-F(v )\geq-U(v-u)$, which along
with \eqref{7.24.1} yields  the Lipschitz continuity of $F$.
 In addition, if $F$ is differentiable at $v\in \bR^{n}$,
then $F(u)-F(v)=D_{u^{i}}F(v)(u^{i}-v^{i})+o(|u-v|)$ 
and using \eqref{7.24.1}
divided by $|u-v|$ and setting $u-v\to0$ we get
$$
  D_{u^{i}}F(v)\omega^{i}\leq \max_{a\in A}[a^{i}\omega^{i}]
$$
for any unit $\omega$. This is only possible if $D  F(v)\in A$
owing to the fact that $A$ is closed, bounded, and convex.
The lemma is proved.


\begin{thebibliography}{mm}

 

\bibitem{Caf89} L.A. Caffarelli, {\em Interior a priori estimates for
solutions of fully non-linear equations\/}, Ann. Math., Vol. 130 
(1989), 189--213.

\bibitem{CC95} L.A. Caffarelli, X. Cabr\'e, ``Fully nonlinear elliptic
equations'',
 American Mathematical Society, Providence, 1995.

\bibitem{CKLS99} M. G. Crandall, M. Kocan, P.L. Lions,
and A. \'Swi{\c e}ch, {\em
Existence results for boundary
problems for uniformly elliptic 
and parabolic fully nonlinear equations\/}, Electron. J.
Differential Equations 1999, No. 24, 1-20,
 http://ejde.math.unt.edu

\bibitem{CKS00} M. G. Crandall, M. Kocan, A. \'Swi{\c e}ch,
{\em $L^p$-theory for fully nonlinear uniformly parabolic equations\/},
 Comm. Partial Differential Equations, Vol. 25  (2000),
 No. 11-12, 1997--2053.



\bibitem{DKL} Hongjie Dong, N.V. Krylov, and Xu Li,
{\em On fully
nonlinear elliptic and parabolic equations in domains with VMO
coefficients\/},  Algebra i Analiz, Vol. 24 (2012), No. 1,
54--95 in Russian; English translation in
St. Petersburg Math. J., Vol. 24 (2013), 39-69.

\bibitem{Es93} L. Escauriaza, {\em $W^{2,n}$ a priori estimates for
solutions to fully non-linear equations\/}, Indiana Univ. Math. J.,
Vol. 42 (1993), No. 2, 413--423.

\bibitem{JS_05} R. Jensen and A. \'Swi{\c e}ch, {\em
Uniqueness and existence of maximal and minimal solutions 
of fully nonlinear elliptic PDE\/}, Comm. on Pure Appl. Analysis,
Vol. 4 (2005), No. 1, 199--207.

\bibitem{Kr85} N. V. Krylov,   Nonlinear elliptic and
 parabolic equations of second order, Nauka, Moscow, 1985 in Russian;
English translation: Reidel, Dordrecht, 1987.

\bibitem{Kr11} N.V. Krylov,
{\em On a representation of fully nonlinear elliptic
operators in terms of pure
second order derivatives and its applications\/},  
 Problemy  Matemat. Analiza,
Vol.
 59, July 2011, p. 3--24 in Russian; English translation:
Journal of
Mathematical Sciences, New York,  Vol. 177 (2011), No. 1, 1-26.

\bibitem{Kr12} N.V. Krylov, {\em
On the existence of smooth 
solutions for fully nonlinear elliptic
equations with measurable ``coefficients''
without convexity assumptions\/},
  Methods and Applications of Analysis,
Vol. 19 (2012), No. 2,   119--146.

\bibitem{Kr13} N.V. Krylov, {\em On the existence of $W^{2}_{p}$ 
solutions for fully nonlinear elliptic
equations under  relaxed convexity assumptions\/}, 
Comm. Partial Differential Equations, Vol. 38 (2013), No. 4,  
687--710.

\bibitem{Kr13.1} N.V. Krylov, {\em
An ersatz existence theorem for fully nonlinear
parabolic equations without convexity assumptions\/}, 
SIAM J. Math. Anal., Vol. 45 (2013), No. 6, 3331--3359. 

\bibitem{Kr12.1}  N.V. Krylov,
{\em Rate of convergence
 of difference approximations 
for uniformly nondegenerate   elliptic Bellman's
equations\/},  
Appl. Math. Optim.,  
Vol. 69 (2014), No. 3,   431--458.

\bibitem{Kr_14}  N.V. Krylov, {\em On $C^{1+\alpha}$ 
regularity of solutions of Isaacs parabolic equations
with VMO coefficients\/}, 
 Nonlinear Differential Equations and
Applications, NoDEA, Vol. 21 (2014), No. 1,  63--85.

\bibitem{Kr_14_1}  N.V. Krylov, {\em To the theory of viscosity solutions
for uniformly elliptic  Isaacs equations\/},
 Journal of Functional Analysis, Vol. 267 (2014), 
 4321--4340

\bibitem{Kr_15} N.V. Krylov, {\em Approximating the value functions
for stochastic differential games   with the ones
having bounded second derivatives\/},
 Stoch. Proc. Appl. Vol. 125 (2015), No. 1,  254--271.

\bibitem{Kr_15_1} N.V. Krylov, {\em To the theory of viscosity solutions for
uniformly parabolic Isaacs equations\/},
 Methods and Applications of Analysis,
Vol. 22 (2015), No. 3, 259--280.

 \bibitem{KT} H.-J,  Kuo and N.S. Trudinger, {\em
Discrete methods for fully nonlinear elliptic equations\/}.
  SIAM Journal on Numerical Analysis, 
Vol. 29 (1992), No. 1,  123--135.

\bibitem{Wa92} L. Wang, {\em On the regularity of fully nonlinear
parabolic equations: I\/},  Comm. Pure Appl. Math., Vol. 45 (1992),
27--76.

\bibitem{Wi09} N. Winter, {\em $W^{2,p}$ and $W^{1,p}$-estimates at the
boundary for solutions of fully nonlinear, uniformly elliptic equations\/},
 Z. Anal. Anwend., Vol. 28  (2009), No. 2, 129--164.

\end{thebibliography}
\end{document}